\documentclass[11pt]{article}
\usepackage{fullpage,subeqn}
\usepackage{epsfig}

\begin{document}

\newcommand{\eq}{\begin{equation}}  
\newcommand{\en}{\end{equation}}
\newcommand{\eqar}{\begin{eqnarray}}
\newcommand{\enar}{\end{eqnarray}}
\newcommand{\ba}{\begin{array}}
\newcommand{\ea}{\end{array}}
\newcommand{\bds}{\begin{description}}
\newcommand{\eds}{\end{description}}
\def\be{{\bf e}}
\def\bn{{\bf n}}
\def\bq{{\bf q}}
\def\bR{{\bf R}}
\def\bx{{\bf x}}
\def\cC{{\cal C}}
\def\cD{{\cal D}}
\def\cL{{\cal L}}
\def\cP{{\cal P}}
\def\cR{{\cal R}}
\def\cS{{\cal S}}
\def\nn{\nonumber}
\def\qand{\quad \hbox{and} \quad}
\def\qie{\quad \hbox{i. e.} \quad}
\def\qfor{\quad \hbox{for} \quad}
\def\qor{\quad \hbox{or} \quad}
\def\qwith{\quad \hbox{with} \quad}

\title{Critical Solutions of Three  Vortex Motion in
the Parabolic Case} 
 
\author{
{\large L. Ting}\\ 
New York University, Courant Institute of 
Mathematical Sciences\\ New York, NY 10012, USA
\and
{\large O. Knio}\\
The Johns Hopkins University, Department of 
Mechanical Engineering\\ Baltimore, MD 21218, USA 
\and
{\large D. Blackmore}\\ 
New Jersey Institute of Technology, Department of 
Mathematical Sciences\\ Newark, NJ 07102, USA} 
 
\maketitle

\footnotetext{Correspondence, 
\textit{E-mail address}: ting@cims.nyu.edu}

\begin{abstract} 

\noindent 
Gr\"{o}bli (1877)   laid the foundation for
the analysis of the  motion of three point vortices 
in a plane by deriving governing equations for triangular configuration
of the vortices. Synge (1949) took this formulation one step further 
to that of 
a similar triangle of unit perimeter, via trilinear coordinates.
The final reduced problem is
governed by  an integrable  two-dimensional system of 
differential equations 
with  solutions represented as  planar trajectories.
Another  key to Synge's analysis was his classification of 
the problem into three
distinct cases: elliptic, hyperbolic and parabolic corresponding,
respectively, to the sum of products of pairs of vortex strengths 
being positive, negative or zero. 
The reduction of the vortex configuration, a curve in space
to a planar curve is one-to-one,
except along a critical planar curve $\cal C$ in the parabolic case.
Each point on $\cal C$ represents  a triangle of unit perimeter 
corresponding to a family of similar vortex configurations, 
expanding or contracting. The latter would lead to coelescence
of the three vortices.
Tavantzis and Ting  (1988) filled most of the gaps
left by Synge regarding the dynamics of the
problem, and  showed in particular 
that points on $\cal C$  corresponding
to   similar expanding families of vortex configurations
are stable while those
corresponding to  similar contracting families are unstable.
Their investigations yielded an exhaustive description of the 
motion and stability of three vortices in a plane  except for the global 
behavior of the vortex configurations 
in a narrow strip containing $\cal C$.
The main contribution of this paper is  a complete description of 
the global dynamics in such a strip,
which  emphaticallly demonstrates that 
three distinct vortices almost never coalesce.
\end{abstract} 

\medskip 
 

\bigskip

\section{Introduction} 
\label{sec-intro}

\begin{figure}[hbt]
\centerline{
\includegraphics[angle=0,width=4.5in, draft=false]{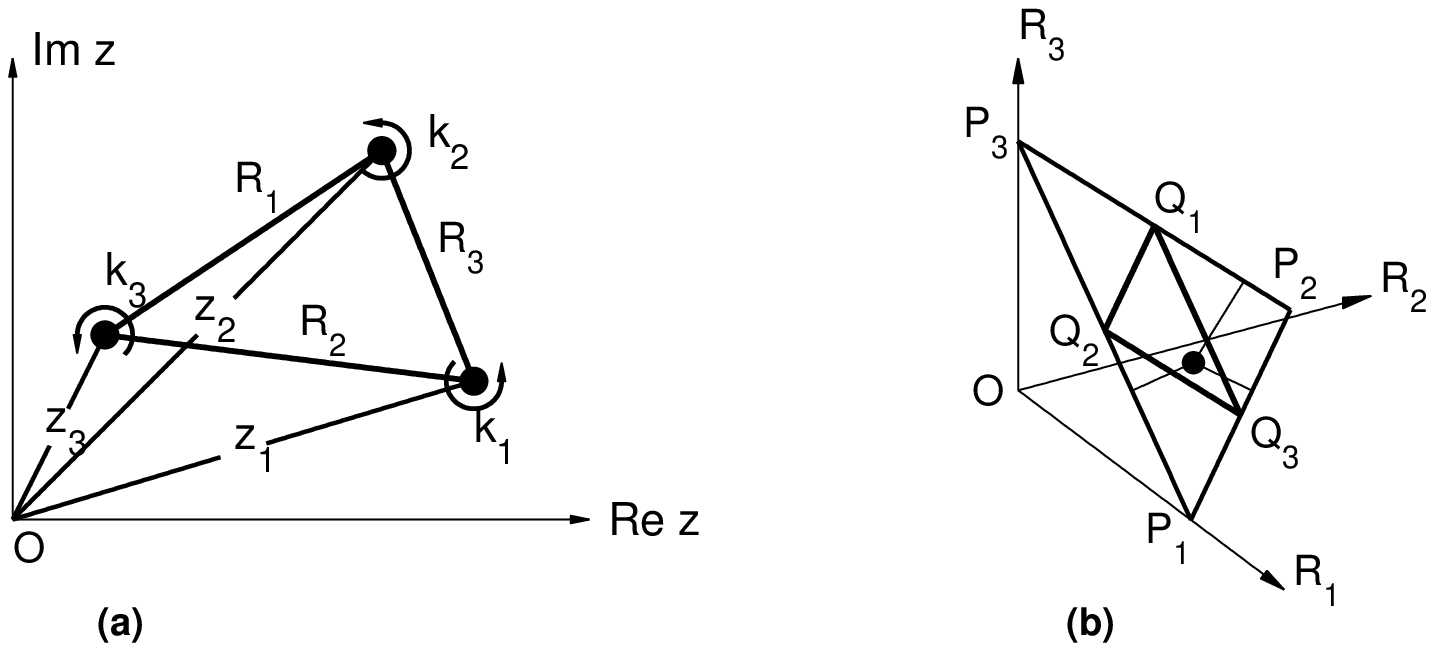} }
\caption{(a)
The triangle $\triangle_p$ formed by
the three point vortices at $z_j (t), \ j=1, 2, 3$,
in the complex $z$-plane, 
and (b) the coordinate axes, $R_1, R_2$ and $R_3$, the plane $\cP$ or 
the $\alpha \beta$ plane
and the  trilinear coordinates  $x_j$, 
shown as the lines from a point, the heavy dot, on $\cP$  to its edges
in the first octant.} 
\label{fig:three}
\end{figure}

The planar motion of three point vortices in an incompressible 
fluid was  studied by 
Gr\"{o}bli in 1877 \cite{grobli} and Synge in 1949 
\cite{synge} among others.
With the aid of the integral invariants, 
Gr\"{o}bli  uncoupled the autonomous six-degree of freedom 
problem to the study the vortex configuration
and the translation and rotation of the configuration and essentially
demonstrated the integrability of the system.
The vortex configuration problem refers to the deformations
of the triangle, $\triangle_p (t)$, formed by the three 
point vortices at 
$z_j(t)$ in the complex $z$-plane. 
The triangle is defined by the lengths,
$R_j (t), \ j=1, 2, 3$, of its sides,  with the  $j$-th side 
facing the $j$-th vortex.
The subscript $p$ denotes the perimeter, $p = R_1 + R_2
+  R_3$, shown in  Fig.~\ref{fig:three} (a).
It is a three degree of freedom problem and 
is integrable, with  the integral curves, or spatial
trajectories, $\bR (t)$,
in  Cartesian coordinates,  $R_j(t)$, given by Synge
\cite{synge}.
The solution was  further symplified by Synge \cite{synge} who 
introduced the trilinear coordinates $x_j , \ j = 1, 2, 3$
with $x_j = R_j /p $. As shown in Fig.~\ref{fig:three} (b), 
the  spatial  integral curve, $\bR (t)$ is projected radially  
onto the  plane, $\cP$, 
which intersects the $j$-th axis at $P_j$ with intercept $\sqrt {2/3}$.
Its cross-section  in the first octant  is the
equilateral $\triangle P_1P_2P_3$ with side $2/ \sqrt 3$
and height $1$, as shown in Fig.~\ref{fig:three} (b).
The heavy dot on the plane $\cP$
with position vector $\sqrt{2/3} \bR /p $,
denotes the radial projection of  $\bR $. The  trilinear
coordinates, $x_j$'s, in turn denote
the distances from the heavy dot  to the sides of the
triangle $P_1P_2P_3$ opposite to  the vertices $P_j$'s,
as shown by the thin lines.
Note that the sum, $x_1 + x_2 + x_3 = 1$, 
is the height of the triangle $P_1P_2P_3$. The 
trilinear coordinates  represent  the sides of  
$\triangle_1 (t)$ with perimeter 1, which is 
similar to and has the same orientation 
as $\triangle_p$ in Fig.~\ref{fig:three} (a).

Because of the triangle inequality, the admissible solutions, $x_j$
of $\triangle_1$, are confined to  $\triangle Q_1Q_2Q_3$
and a point on its positive (negative) side with normal vector 
pointing away from (towards)  the origin $O$ corresponds
to a configuration
$\triangle_p$
with  the three vortices in counterclockwise (clockwise)
orientation. Two points on the opposite sides of $\cP$
with the same $x_j$ 
are called the image points of each other.

Synge \cite{synge} identified
the sum of the products of the vortex strengths, $k_j$, 
\eq
K = k_1 k_2 + k_2 k_3 + k_3 k_1,  
\label{eq:K}
\en 
as the primary parameter, and described  the topologies
of the integral curves and critical points in the  plane $\cP$
according to the classifications,
\eq
\hbox{elliptic;}\ K > 0; \qquad  \hbox{parabolic;} \ K = 0;
\qand  \ \ \hbox{hyperbolic}, \ K < 0.
\label{eq:class}
\en
In particular, he studied 
the critical point, {\bf E}, at the centroid of
$\triangle P_1P_2P_3$ and also that of $\triangle Q_1Q_2Q_3$
($\triangle_Q$ for short hereafter) on
the positive side of $\cP$,
where $ x_j = 1/3$. The corresponding 
configuration  is an equilateral triangle and was shown to be
stable (unstable) for  the elliptical   (hyperbolic) cases.
This is true also for the image of  $E$, the centroid $E^*$
on the negative side of $\cP$.

For the parabolic case, Synge pointed out that
the critical point $E$ lies on a critical curve $\cC$ and
partitions $\cC$ into expanding and contracting branches,
while the images of those two branches on $\cC^*$ are
contracting and expanding, respectively.
The study of the local stability and
global trajectories near these two branches on $\cC$  and their images
is  the main subject of this paper. 

The planar vortex problem was  studied  
via Hamiltonian formalism  by  Lin in 1943, \cite{lin43},
Novikov in 1975 \cite{novikov}, Aref in 1979 \cite{aref} among others. 
Tavantzis and Ting in 1988  \cite{tavtin} continued  
Synge's analysis via trilinear coordinates.
For each classification by $K$ in (\ref{eq:class}), they
located all the critical points and   studied their 
stability  and described
the topologies of the integral curves 
and  separatrices in   the plane $\cP$.

Numerical studies of the motions of three  vortices 
in a half-plane were carried out by Knio et al. \cite{knicol}. 
for  vortices of equal strength, $k$, i.~e., $K = 3k^2 > 0$ in
the elliptic cases. The numerical results showed that
the motion of the vortices
can be regular or chaotic depending on  their
initial configuration, which was taken to be either an equilateral
triangle or collinear. 
The transition  from regular to chaotic
for those  half plane  problems can be anticipated when
the  distances between the vortices  are much  smaller than 
their distances from the boundary. The half plane problems 
can then be  considered
as perturbations
of three  vortices in the entire plane, 
for which the equilateral
configuration corresponds to a center while the collinear
configuration  corresponds to a saddle
point in elliptic cases, \cite{synge}. 
This observation motivated  a sequence of theoretical and numerical
studies of perturbed three vortex problems in parallel with
the classical perturbed two body problems. Here the perturbations
arise  for example
from the presence of  a boundary and/or additional vortices far away
from the  three main vortices
or from  the approximations the interactions of 
three coaxial vortex rings by three  vortices in a meridian plane
when the distances between the rings are much smaller than 
the radii of the rings. See the  article
by Blackmore et al. \cite{BTK06} and the references therein 
and recent studies
by Ting et al. \cite{TKB07} and Knio et al. \cite{KBT08}.
Several numerical examples were presented to show transitions 
from periodic and quasi-periodic to chaotic regimes  
in accordance with the theoretical results.

This paper does not deal with perturbations. Rather, it presents an
exhaustive study of the global behavior of the 
dynamics of three  vortices in a plane in  parabolic cases,
with initial vortex configurations in a narrow strip containing
the critical curve $\cC$. This completes  the local stability
analysis of $\cC$ presented Tavantzis and Ting
\cite{tavtin}. 

For the strengths of three  vortices,  at least 
two of them  have to have the same sign. Thus
we can always choose the sign of those two to be positive and denote 
them by $k_1$ and $k_2$ with 
\eq
k_1 \ge k_2 > 0, \qand   k_3= - k_1k_2/(k_1+k_2) < 0,
\label{eq:order}
\en
for the parabolic case, $K=0$. 
After ordering  the three  vortices by Eq.~(\ref{eq:order}), 
the initial positions of 
the vortices $z_j(0)$ in the $z$-plane define the sides $R_j(0)$
of the triangle $\triangle_p$ and hence the trilinear coordinates $x(0)$
in the plane $\cP$ and in addition 
the orientations of the vortex configuration, 
clockwise or counterclockwise,
denoted by the index $\gamma = +1$ and $-1$, respectively.
To account for the opposite orientations, $\gamma = \pm 1$,
Synge made use of both sides, the positive and negative sides,
of $\cP$, or  $\triangle P_1P_2P_3$. The normal vector to the positive
or the negative side, points away from or towards the origin $O$, 
respectively, (see Fig.~\ref{fig:three}).
He assigned  the points $x_j$ with $\gamma = +1\ (\gamma = -1)$ to be 
on the positive (negative) side of 
$\triangle P_1P_2P_3$. A change of orientation can take place only
when crossing over an  edge of the triangle.

It was pointed out by Synge \cite{synge} that the
mapping from a point $\bf R$ on the spatial integral curve
to its radial projection  onto a point on the
planar  curve in  $\cP$, 
is one to one if the velocity $\dot {\bR} (t)$ is not in the 
radial direction, i.~e., $\dot {\bR} \times \bR \not= 0$,
which can occur
only on a planar  critical curve $\cC$ in a parabolic case.
The critical  curve is partitioned into  expanding and contracting branches.
Corresponding to a point  on  the curve, $\cC$,
the sides $x_j$'s of the $\triangle_1$ remain stationary,
while $\triangle_p(t)$, similar to $\triangle_1$, is expanding, with
$\dot p > 0$ on one branch and  contracting $\dot p < 0$  on the adjacent
branch.
It was shown in \cite{tavtin} that the similar expanding
solution is stable while the contracting one is unstable.
Similar contracting solutions, albeit unstable,  have 
received  attention
because they could in principle lead to simultaneous colliding or
coalescence of three  vortices (see for example Newton  
\cite{newton}, pp. 84-88).
Thus we need to reinforce the local instability analysis  of 
\cite{tavtin}, for a point on a contracting branch of $\cC$
by nonlocal analysis leading to the conclusion that
coalescence of three vortices is unlikely.

In Sec.~\ref{sec-background} we present a
brief review of the  three
point vortex problem, define the symbols, and summarize the relevant
results of
Synge \cite{synge} and Tavantzis and Ting \cite{tavtin}.
To render this paper self-contained, we reproduce derivations 
and descriptions of
results relevant to the current investigation,
namely, the spatial trajectories of vortex configurations
and their reduction to planar trajectories  via trilinear coordinates.

In Sec.~\ref{sec-traj-cC}, we study in detail the planar 
trajectories and the critical curve $\cC$ in parabolic cases.
We extend the local stability analysis of the  critical curve
in \cite{tavtin}
by showing that a small  deviation off $\cC$ from a point 
on a contracting branch produces a trajectory that departs
from the neighborhood of $\bx$ and is  finally
attracted  to an expanding branch of the critical curve.
Depending on the location of the point $\bx$ on $\cC$ and the 
perturbation of $\bx$ to a point above or below $\cC$, we
identify three types of trajectories. An outline of 
this investigation was presented recently
by Ting and Blackmore \cite{tinbla}.

In Sec.~\ref{sec-numerical} we present numerical 
examples to demonstrate how the 
the configuration $\triangle_1$ 
departs from the neighborhood of a contracting branch of
the critical curve $\cC$  along  a trajectory 
to the expanding branch. We show the three types
of trajectoties and the variations of the configurations
depending on the locations of
initial point $\bx$, as predicted
in Sec.~\ref{sec-traj-cC}.

\section{Background}
\label{sec-background}
Adopting the symbols of \cite{synge}  and \cite{tavtin}, 
we relate the strength of a
vortex  $k_j$ to its circulation $\Gamma_j$   by $ k_j = \Gamma_j/2
\pi, \ \ j = 1, 2, 3 $, and denote the position of the $j$-th
vortex  in the
complex plane by   $ z_j(t) = \Re z_j + i \Im z_j$,
and the length of the  side of $\triangle_p$ facing its $j$-th
vertex at $z_j$ by $R_j $. 

Without loss of generality, we  assign
the strengths of the three vortices according to  Eq.~(\ref{eq:order}).
We set the initial perimeter as the length scale, $p(0) = 1$ and then 
the time scale is defined by $p^2(0)/k_2$.
The equations of motion of three point vortices 
in a planar incompressible inviscid flow are:
\eq
 \dot z_j(t) = -i \sum_{m \not= j} \frac{k_m}{\bar z_m - \bar z_j},
 \qfor j= 1, 2, 3.
\label{eq:motion}
\en
The summation in $m$ ranges over $1, 2$, and 3.
We obtain from Eq.~(\ref{eq:motion})
the well-known theory of Kirchhoff or the conservation laws of
vorticity and the moment of vorticity
in the $z$-plane (\cite{lamb}, pp. 229-230),
\eq
\sum_j k_j \dot z_j = 0 ,\qquad \hbox{and} \qquad
\sum_j k_j \bar z_j \dot z_j(t) = i K
\label{eq:conserv1-2}
\en
The second equation in polar coordinates yields
\eq
\sum_j k_j \rho_j^2  = \hbox{const.} \quad\hbox{and}\quad
 \sum_j  k_j \rho_j^2 \dot \theta_j =   K,
\label{eq:polar}
\en
where $\rho_j$ and $\theta_j$ denote the magnitude and argument of $z_j$,
respectively.

In the next section, we will reproduce the analyses of Gr\"{o}bli (1877)
and Synge (1949) which decoupled 
the six equations for  the real and imaginary parts of 
Eq.~(\ref{eq:motion}), $j=1, 2, 3$, to three equations
for $R_j$, the vortex configuration and three for the motion
of the configuration. The equations for the latter will come
from the above two conservation equations in  (\ref{eq:conserv1-2}).

\subsection{Deformation of the configuration $\triangle_p$}
\label{sec-R1R2R3}

The deformations of $R_j$'s  are governed by  linear
combinations of Eq.~(\ref{eq:motion}).
For example, from a linear combination of the first two equations  
(\ref{eq:motion}), we get 
\eq
(\bar z_1 - \bar z_2) (\dot z_1 - \dot z_2)   =  R_3 \dot R_3 + 
i  R_3^2 \dot \psi_{21} =   i (k_1 + k_2) -i  k_3 R_3 
[R^{-1}_2 e^{i [\psi_{13} - \psi_{21}]}  
-R^{-1}_1 e^{i [\psi_{23} - \psi_{21}]}] .
\label{eq:dz1z2}
\en
Here   $\psi_{jl}$ denotes the  argument of the side $z_l - z_j$, and
$\psi_{lj} =  \psi_{jl} + \pi$ denotes that of $z_j - z_l$.
The real part of the equation 
yields the rate of change of the side $R_3 = |z_1 -z_2|$.
Likewise, Synge \cite{synge}  obtained
the rates of change of all three sides of $\triangle_p$,
\eq
\frac{\dot R_1}{ k_1 R_1 (R_3^2 - R_2^2)} = 
\frac{\dot R_2}{ k_2 R_2 (R_1^2 - R_3^2)}  = 
\frac{\dot R_3}{ k_3 R_3 (R_2^2 - R_1^2)} =
\frac{2\gamma |A|  }{R_1^2 R_2^2 R_3^2}.
\label{eq:dotRj}
\en 
where $ |A| = [s(s-R_1)(s-R_2)(s-R_3)]^{1/2} $, with $s= p/2$,
stands for the area of $\triangle_p$. 
From the initial data for the $z_j$'s, we have the data for the $R_j$'s  of the
$\triangle_l$  and then
Eqs.~(\ref{eq:dotRj}) define the spatial
integral curve, $\bR(t)$, in the first octant, where $R_j > 0, \ j =1, 2, 3$.
The integral curve $\cR$ is  defined directly by the
two integral invariants,  or by the intersection of two surfaces,
\eqar
k_1^{-1} R_1^2 + k_2^{-1} R_2^2 + k_3^{-1} R_3^2 & = & 
\hbox{initial value}\ a \label{eq:synge1} \\
R_1^{1/k_1}  R_2^{1/k_2}  R_3^{1/k_3}  & = & 
\hbox{initial value} \ b.
\label{eq:synge2}
\enar
They were obtained by Kirchhoff for $N$-vortices, with $N=3$ here,
(see e.~g. \cite[page 230]{lamb}). From  $\dot \bR (t)$, we define
the direction along the integral curve  and call the directed integral
curve the trajectory of the configuration, or simply the trajectory.

To recover the positions $z_j(t)$ corresponding to a point on the 
integral curve $\cR$, we need to find three combinations of
the primary six equations in (\ref{eq:motion})  independent 
of the three equations for the configuration (\ref{eq:dotRj}).
These three combinations appear readily from the real and imaginary
parts of the conservation of the center of 
vorticity in (\ref{eq:conserv1-2}) and
the imaginary part  of the polar moment of vorticity in
Eqs.~(\ref{eq:polar}). The real part of the equation for the polar
moment is not independent because it is 
equivalent to Eq.~(\ref{eq:synge1}), see \cite{synge}.

We first find
the dependence of $R_j$ and their inclinations, say $R_3$ and $\psi_{12}$,
on $t$ by integrating   Eq.~(\ref{eq:dotRj}). That
would then define the temporal variations of the configuration 
$\triangle_p$ and its inclination.
Finally the locations $z_j(t)$ of the vertices of  $\triangle_p(t)$
are defined by the first invariant in (\ref{eq:conserv1-2}), 
or by locating
the (stationary) weighted center of $\triangle_p(t)$, 
$z_c=(k_1z_1(0)+k_2z_2(0)+ k_3 z_3(0))/ (k_1 + k_2 + k_3)$, for $k_1
+k_2 + k_3 \not= 0$, which holds for parabolic cases.

The  reduction   
of the spatial curve $\bR$ to a planar curve by Synge \cite{synge}
via trilinear
coordinates is outlined in the next subsection.

\subsection{The trilinear coordinates}
\label{sec-x1x2x3}

Synge [2] projected the point $\bR$ 
radially onto the point, $\bq$  on the plane $\cP$,
\eq
R_1 + R_2 + R_3 = \sqrt {2/3},
\label{eq:calP}
\en
which intersects  the axes  respectively at $P_1$, $P_2$ and $P_3$, 
with equal intercept, $\sqrt{2/3}$. The plane $\cP$
intersects the coordinate planes at the same 
angle, $\arcsin  \sqrt {2/3}$.
As shown in Fig.~\ref{fig:three}~(b),
the equilateral $\triangle P_1P_2P_3$ denotes the 
section of $\cP$ in the first
octant with each side $2/\sqrt 3$ and height $H =1$. 
The trilinear coordinates $x_j, \ j=1, 2, 3$ of 
point $\bq$ in $\cP$ are
the signed  distances from the sides of $\triangle P_1P_2P_3$. 
The sign of
the distance from a side is 
positive (negative) if the distance from the side to  $\bq$ points
inwards to (outward from) $\triangle P_1P_2P_3$.
The sum of the trilinear coordinates is always equal to the height $H$,
\eq
x_1 + x_2 + x_3 = 1.
\label{eq:tri-x}
\en
The $x_j$'s  are all
positive when $\bq$ lies inside  $\triangle P_1P_2P_3$ and
are related to the Cartesian coordinates $R'_j$ of $\bq$ by, 
\eq
\frac{R'_j}{x_j} =  \sqrt {\frac{2}{3}}, \qquad
\hbox{while} \quad \frac{R_j}{R'_j} = \frac{p}{\sqrt {2/3}}
\quad \hbox{and} \ \ \frac{R_j}{x_j} = p, \quad j = 1,2,3.
\label{eq:similar}
\en
Thus the $x_j$'s
also denote  the sides of  $\triangle_1$ with perimeter $ 1$, 
similar to $\triangle_p$ with  the same orientation.
The mapping of the spatial trajectory $\bR(t)$ to the planar trajectory
$\bx (t) $  is
one-to-one, provided that the 
radial projection of  $\dot \bR$ onto $\cP$ is nonzero,
or the spatial trajectory is not radial 
\eq
\dot \bR \times \bR \neq 0.
\label{eq:one2one}
\en
This condition is fulfilled for the elliptic and hyperbolic cases.
The condition can be violated only for the parabolic case, where 
the spatial trajectory can be radial. 
Then the radial projection of the trajectory 
onto $\cP$ reduces to one singular point $\bq$. 
The loci of these singular 
points  is called the  critical 
curve $\cC$  such that each point $\bq$ on $\cC$ represents
a $\triangle_1$ corresponding to 
a radial trajectory of $\triangle_P \sim \triangle_1$ either moving 
away from  the origin with $\dot p > 0$ or moving inward with $\dot p< 0$,
see \cite{synge} and \cite{tavtin}. 
More details on the  critical curve
will be elaborated later in {\bf Sec.}~\ref{sec-traj-cC}, to set the stage
for the main objective of the current paper, which
is to identify branches of the critical curve that
attract or repel   nearby trajectories.

Due to the triangle inequality,
\eq
x_1 + x_2 \ge \frac{1}{2} \ge x_3, \quad
x_2 + x_3 \ge \frac{1}{2} \ge x_1, \quad
x_3 + x_1 \ge \frac{1}{2} \ge x_2,
\label{eq:inequality}
\en
$\bx$ has to lie in the $\triangle_Q$
with vertices $Q_j,  \ j= 1, 2, 3$ 
lying  at the midpoint of the $j$-th side of $\triangle P_1P_2P_3$.
The equality signs in (\ref{eq:inequality})
hold   on the edges $Q_1Q_2$ , $Q_2Q_3$ and $Q_3Q_1$ of $\triangle_Q$, 
where $x_3,\ x_1$ and $x_2 = 1/2$, respectively,
and the three vortices
are collinear with the area of $\triangle_p = 0$.
The vertex $Q_j$ of $\triangle_Q$   corresponds to a singular case
where the three vortices degenerate to two. For example, 
at $Q_1$, where $x_1=0$, and  $ x_2 = x_3 = 1/2$, 
or $R_1=0$ and $R_2=R_3$,
the three vortices degenerate to two, one at $z_1$ with strength $k_1$
and one at $z_2=z_3$ with strength $k_2 + k_3$.

On the positive (negative)  face of the plane $\cP$, or the $\triangle_Q$, 
where the  normal vector to the face  points away  from  (towards) the origin, 
we  assign $\gamma=1$ ($\gamma = -1$).
With the same trilinear coordinates $(x_1, x_2, x_3)$, denoted by $\bx$,
the point on the positive (negative) face of $\triangle_Q$
implies  that the vertices of $\triangle_1$, or the three vortex centers
$z_1, z_2$ and $z_3$
in the $z$-plane, are in the counterclockwise (clockwise) orientation,
Fig.~\ref{fig:three} (a) shows the vortex centers with
$\gamma = +1$. We call the point $\bx$ on the negative side of $\cP$
the image of $\bx$ on the positive side and vice versa.
The orientations of the vortex centers or $\triangle_1$
can change only in crossing over an edge of $\triangle_Q$,
where $\triangle_1$ is collinear.

In trilinear coordinates, the two invariants (\ref{eq:synge1}) 
and (\ref{eq:synge2})  become,
\eqar
& & [k_1^{-1} x_1^2 + k_2^{-1} x_2^2 + k_3^{-1} x_3^2] p^2 =  a \qand 
\label{eq:synge3a} \\
& & x_1^{1/k_1}  x_2^{1/k_2}  x_3^{1/k_3}  =  b p^{K/(k_1k_2k_3)}.
\label{eq:synge3b}
\enar
With the  elimination of  the perimeter $p$, these two equations yield 
an integral curve, or a trajectory, in $\cP$,
valid for  all $K$. It  is
\eq
[ \frac{x_1^2}{k_1} + \frac{x_2^2}{k_2}+ \frac{x_3^2}{k_3} ]^{K/(2 k_3)}
[x_1^{1/k_1} x_2^{1/k_2}  x_3^{1/k_3} ]^{k_1k_2} = \hbox{const.} \bar I.
\label{eq:tt1}
\en
For a parabolic case, $K=0$,  Eq.~(\ref{eq:tt1}) reduces to
\eq
x_1^{k_2} x_2^{k_1} x_3^{k_1k_2/k_3} = \bar I 
\quad \hbox{or} \quad (\frac{x_1}{x_3})^{k_2} (\frac{x_2}{x_3})^{k_1} = \bar I,
\label{eq:tt2}
\en
which is equivalent to Eq.~(\ref{eq:synge3b}) for the parabolic case when
the constant $\bar I$ is identified as $b^{k_1 k_2}$. 
In the following  subsection we recount the critical points 
in the planar trajectories first presented by Synge \cite{synge}
with special attention to the parabolic cases.

\subsection{Critical points $\bR$ in space}
\label{ssec-critical}

It was shown in \cite{synge} that
at  a critical  point  $\bR$, the configuration $\triangle_p$
is stationary, i.~e., $\dot R_j = 0, \ j =1,
2, 3$,  if and only if   $\triangle_p$ is either equilateral  with
$R_1=R_2=R_3 = p/3$ or has  zero area, $|A| = 0$. The latter
requires that the configuration $\triangle_1$ is collinear lying on an
edge of $\triangle_Q$ or at the vertices of $\triangle_Q$, where 
two vortices coincide to one,
see Eqs.~(\ref{eq:dotRj}). 
These critical points of $\bR$ in space 
are also critical points $\bq$
in plane $\cP$, or $\bx$
in the trilinear coordinates, $x_1, x_2, x_3$ with  $p=1$, 
see Eqs.~(\ref{eq:similar}). 
An equilateral $\triangle_p$, implying $x_1=x_2=x_3 = 1/3$
corresponds to the centroid $E$  of $\triangle_Q$ for positive
orientation , $\gamma = 1$ or 
its image $E^*$ on the opposite side for $\gamma = -1$. 
For an  equilateral configuration, it was shown that   $E$ and  $E^*$,
are stable for the elliptic case 
and unstable for the hyperbolic
case \cite{synge}. For the 
parabolic case, $E$ and $E^*$ are degenerate singular points 
because they lye on the critical curve, which will
be discussed in {\bf Sec.}~\ref{sec-traj-cC}.

\begin{figure}[hbt]
\centerline{
\includegraphics[angle=0,width=3.0in, draft=false]{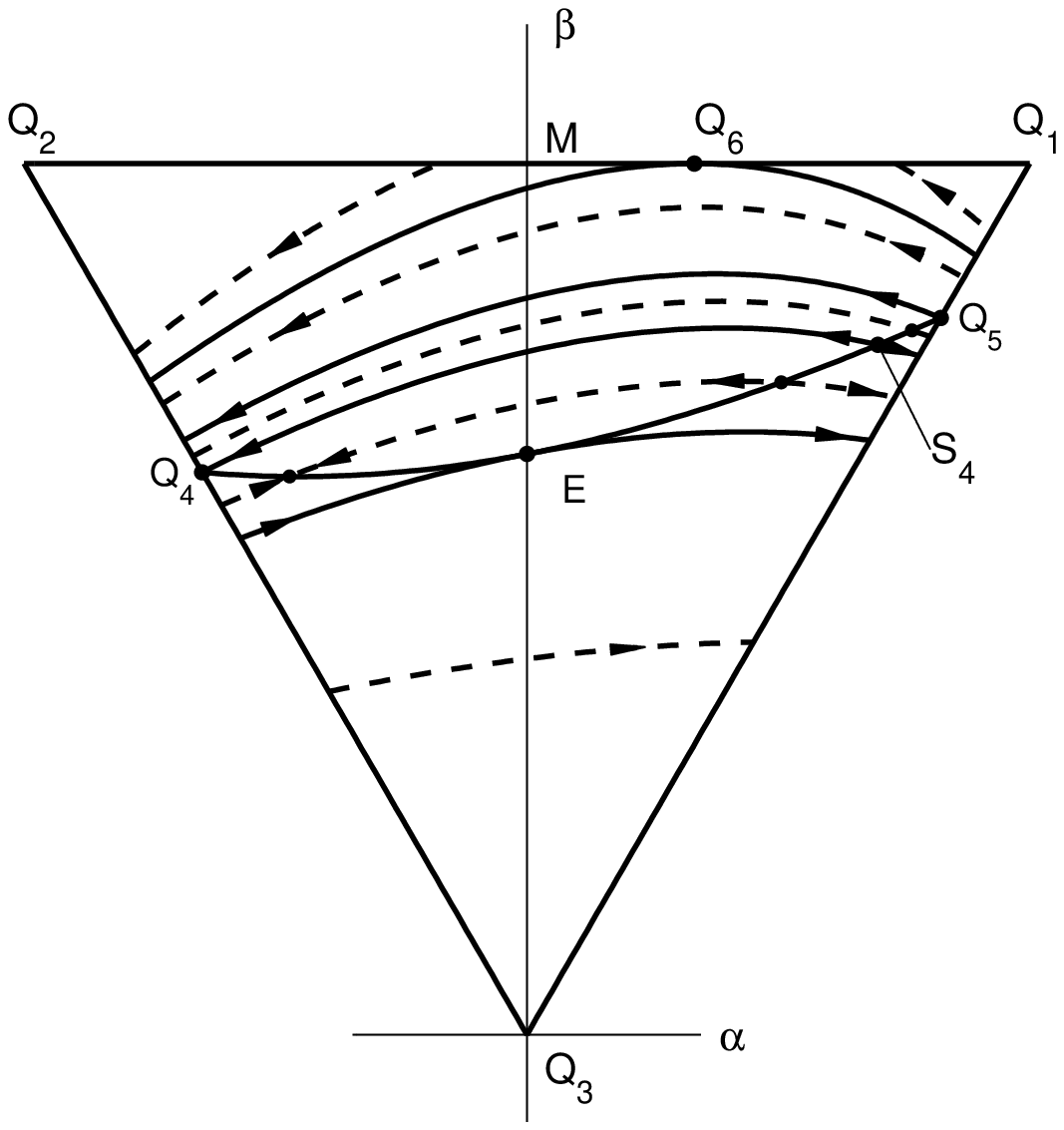} }
\caption{The trajectories for the parabolic case, 
$k_1k_2 + k_2k_3 + k_3 k_1 = 0$,
on  the  positive side of the plane $\cP$ or 
the equilateral $\triangle_Q$.
That is, the configuration
$\triangle_1$ has  counterclockwise orientation. 
Here $k_1=2,\  k_2 =1$
and $k_3 = - 2/3$. Also shown are the  coordinates, $\alpha \beta$,
with the vertex $Q_3$ as the origin  and  $ Q_3 E$
as the $\beta$ axis.}
\label{fig:parabA}
\end{figure}

For a collinear configuration,  the vortices
$z_j$'s have to lie on a straight line, $\cL$. 
Hence we could use real numbers
$s_j$  to locate $z_j$ along  $\cL$, by setting   $s_1= 0$
and assigning its direction such that
$s_2 = s_2 -s_1 = R_3 > 0 $. 
The condition of collinearity  is
\eq
z_3 - z_1 = \nu ( z_2 - z_1) , \qor s_1 = 0, \ \ s_2 = R_3 > 0, \quad 
s_3 = \nu s_2.
\label{eq:collinear}
\en 
where $\nu$ is a real number.
The last equation implies $R_2 = |\nu| R_3$.
From the values  of $\nu$, whether negative, in $ [0, 1]$
or greater than $1$, we see  the relative positions of the vortices
along $\cL$, the collinear configurations
and the locations of the critical points
on the edges of $\triangle_Q$
(see Fig.~\ref{fig:parabA} and \cite{tavtin}).
For Eq.~(\ref{eq:collinear}) to hold for all $t$, we obtain from
the equations of motion (\ref{eq:motion}), that $\nu$ has to be the
root of the cubic equation, $f(\nu) = 0$, where
\eq
f(\nu) = (k_1 + k_2) \nu^3 - (k_1 + 2 k_2) \nu^2
- (k_1 + 2 k_3) \nu + (k_1 + k_3) 
\label{eq:fnu}
\en
It was shown in \cite{tavtin}, that  
there is a value $K_* < 0$, such that
the cubic equation (\ref{eq:fnu}) 
has three  unequal real roots for $K > K_*$ and only one
real root for $K < K_* < 0$ in hyperbolic cases.
For the parabolic cases, we have  explicit formulas for the
three roots of $f(\nu) = 0$.
With $k_3 = -k_1k_2/(k_1+k_2)$, we see that 
 $k_1/(k_1 + k_2)$ is a root of
Eq.~(\ref{eq:fnu})  and hence the three roots are:
\eqar
& & \nu = k_1/(k_1 + k_2) \in [1/2, 1), \quad \ Q_6 \ \hbox{on side}\ Q_1Q_2,
\label{eq:root1}\\
& & \nu = [k_2 - \sqrt{k_1^2 + k_1k_2 + k_2^2}]/(k_1 + k_2) < 0,
\ \ Q_4 \ \hbox{on side}\ Q_2Q_3, \label{eq:root2} \\
& & \nu = [k_2 + \sqrt{k_1^2 + k_1k_2 + k_2^2}]/(k_1 + k_2) > 1,
\ \ Q_5 \ \hbox{on side}\ Q_3Q_1. \label{eq:root3}
\enar
These  three critical points lie on the respective sides
of $\triangle_Q$.

The trajectories defined by Eq.~(\ref{eq:tt1}) and the critical points  
were shown in Figs. 3, 4, and 5  in \cite{tavtin} for typical  
elliptic , hyperbolic and parabolic cases, respectively.  
In this paper we study the critical solutions, which 
occur only in the parabolic case; therefore,
the figure for
the parabolic case in \cite{tavtin}
is reproduced here in Fig.~\ref{fig:parabA}.

It was noted before that a stationary point 
in the $\bR$ space corresponds to
a stationary point $\bq$ on the plane $\cP$,
but the converse may not be true. A  stationary point
$\bq$ on  $\cP$ could 
correspond to similar $\triangle_p$ in the $\bR$ space with perimeter
 $p(t)$ varying. Those stationary points on $\cP$
yields the critical curve $\cC$.
We shall  examine in detail
the geometry of the critical curve and 
the planar  trajectories, in particular those intersecting
the critical curve, in the next section.

\section{The planar trajectories and critical curve for the parabolic case}
\label{sec-traj-cC}

A point $\bq$ on the plane $\cP$, represented by 
trilinear coordinates $x_j, \ j=1,2,3$, is related to the 
coordinates  $R_j, \ j=1, 2, 3$ of the point in space 
by Eq.~(\ref{eq:similar}). We have
$\dot R_j = p \dot x_j + x_j \dot p$, 
and  convert  Eq.~(\ref{eq:dotRj}) for $\dot R_j$ 
to equations for $\dot x_j$ and $\dot p$. They are:
\eqar
& & \dot x_1 /x_1 = -  \dot p / p 
+ \gamma{\cal H} k_1  (x_3^2 - x_2^2), 
\qquad\quad  
\dot x_2 / x_2 = - \dot p / p + \gamma{\cal H} k_2  (x_1^2 - x_3^2) 
\label{eq:dotx1-2} \\
& & \dot x_3 / x_3= -  \dot p / p + \gamma{\cal H} k_3  (x_2^2 - x_1^2), 
\qquad\quad \hbox{where}\ \ {\cal H} = 
2  |A| / (p^2 x_1^2 x_2^2 x_3^2),
\label{eq:dotx3-H} \\
&\hbox{and}\quad&   \dot p/ p  = \gamma {\cal H}[k_1x_1 (x_3^2 - x^2_2)
+ k_2x_2 (x_1^2 - x_3^2) + k_3 x_3 (x_2^2 - x_1^2)].
\label{eq:dotp} 
\enar 
The last equation (\ref{eq:dotp})
follows  either directly from Eq.~(\ref{eq:dotRj}) or from the sum of
the preceding first three equations multiplied by $x_j, j=1,2,3$, 
respectively.

We eliminate $\dot p/p$ from the first three equations and obtain
\eqar
& & \frac{\dot x_1}{x_1} - \frac{\dot x_2}{x_2}  
= - \frac{\gamma}{k_3} {\cal H}(\bx) {\cal Y}(\bx) 
\qquad\hbox{and}\quad
\quad \frac{\dot x_1}{x_1} - \frac{\dot x_3}{x_3} 
= \frac{\gamma}{k_2}  {\cal H} (\bx) {\cal Y}(\bx),
\label{eq:x1x2x3}\\
\hbox{where}\qquad& & {\cal Y} (\bx)
= k_2k_3 x_1^2 + k_3k_1 x_2^2 + k_1k_2 x_3^2.
\label{eq:calY}
\enar
Using the identity, $\sum_j x_j = 1$, we get 
$\sum_j \dot x_j = 0$, and  together with
the two equations in (\ref{eq:x1x2x3}) we  solve for $\dot x_1$.
We then get 
\eqar
& & \quad \dot x_j =  \gamma F_j\ {\cal H}\ {\cal Y} \label{eq:FHY} \\
&\hbox{where}& \quad
F_1 = x_1[ \frac{x_3}{k_2}  - \frac{x_2}{k_3}], \ \
F_2 = x_2[ \frac{x_1}{k_3}  - \frac{x_3}{k_1}], \ \
F_3 = x_3[ \frac{x_2}{k_1}  - \frac{x_1}{k_2}].
\label{eq:dotxj}
\enar
To arrive at  the same factor ${\cal Y} (\bx)$ on the right hand sides
of Eqs.~(\ref{eq:x1x2x3})   and 
(\ref{eq:dotxj}), we  used   $K=0$.
Thus, the two linearly independent equations in (\ref{eq:dotxj}) define
the planar trajectory in $\cP$ for parabolic cases only.
These two equations plus
Eq.~(\ref{eq:dotp}) for the variation of the perimeter are
equivalent to the three equations (\ref{eq:dotRj})
for the spatial trajectory in  parabolic cases.

For  points on $\cC$ other than the critical points
$Q_4, E, Q_5$ and $E^*$ we have the common factor ${\cal H} > 0$.
For all the points on $\cC$ to be stationary, 
it is necessary that the other common factor ${\cal Y} =0$,
\eqar
& & {\cal Y}(\bx) = k_2k_3 x_1^2 + k_1k_3 x_2^2 + k_1 k_2 x^2_3 = 0,
\label{eq:cC0}\\
&\hbox{or}\qquad\qquad&  Y (\bx) = -{\cal Y}/k_3 = - k_2 x_1^2 - k_1 x_2^2
+ (k_1 + k_2) x_3^2 = 0.
\label{eq:curve3}
\enar
This is the equation for the critical curve $\cC$,
which exists only in the parabolic case.

Here the  derivation of the equation for $\cC$ follows directly
from the equations for $\dot x_j$ in (\ref{eq:dotxj}) and
differs from the original derivation of Synge \cite{synge}.
Besides being   an alternative derivation,
the equations in (\ref{eq:dotxj})
show explicitly the dependence of $\dot x_j$'s on
$\gamma$ and $\cal Y$ or $Y$. That is {\it the direction of a trajectory 
is reversed when it crosses over an edge of $\triangle_Q$ to the opposite
side as $\gamma$ changes from $\pm 1$ to $\mp 1$, respectively,
or when it crosses over the critical curve $\cC$, where  $Y=0$,
say from above to below as $Y$ changes from positve to negative.}
These results  will be  needed for   
our study of the trajectories in the strip
containing the critical curve in {\bf Subsecs.~\ref{ssec-direc}}
and {\bf \ref{ssec-bifurcation}}.

On the stationary critical curve $\cC$,
the variation of the perimeter  $\dot p$ can be obtained 
directly from one of the equations in (\ref{eq:dotx1-2}), say
the first one. The result is 
\eqar
& & p^2(t) = p^2(0) + 4\gamma D_0 S_0 t,
\label{eq:soln-p}\\
&\hbox{where}\quad&                                         
D_0 = \frac{x_2^2 - x_1^2}{k_1k_2}
\quad \hbox{and} \ \ 
S_0 = 
\frac{\sqrt{(1-2x_1)(1-2x_2)(1-2x_3)}}{2k_1k_2k_3 (x_1x_2x_3)^2}
\quad \hbox{are constants.}\ \ \
\label{eq:S0-D0}
\enar

We note that when $p(t) = p(0)$, $\bR$ 
is  stationary. This  can happen only
at the centroid $E$ and its image $E^*$ where $x_1=x_2=x_3 = 1/3$
and $D_0 = 0$
and at the critical points $Q_4$ and $Q_5$ on the two sides of
$\triangle_Q$ where $|A| = 0$ and hence $S_0=0$.
Thus the critical curve $\cC$ has to pass through $E$ and $E^*$ and
cross over the sides
$Q_2Q_3$ and $Q_3Q_1$ at the critical points $Q_4$ and $Q_5$, respectively
(see Eqs.~(\ref{eq:root2}) and (\ref{eq:root3})).

To describe the planar trajectories  and the critical curve
$\cC$ in the plane $\cP$, i.~e., to find their slopes and concavities,
we need to replace the
trilinear coordinates satisfying
the constraint, $x_1+x_2+x_3 = 1$, by 
Cartesian  coordinates in $\cP$. We choose the coordinates
$\alpha$ and $\beta$ with  the vertex $Q_3$ as the  origin 
and the bisector of
$\angle Q_3$, $Q_3EM$, as the $\beta$ axis, where  $M$ denotes
the midpoint of the side $Q_1Q_2$ (see Fig.~\ref{fig:parabA}), with
\eqar
\beta & = & x_3   \qquad  \alpha = (x_2 - x_1)/ \sqrt 3 ,
\qquad \hbox{with  inverse transformation,} 
\label{eq:xjtoab}\\
x_3 & = & \beta, \ \  x_1 = [1 - \beta - \alpha \sqrt 3] /2, \qand  x_2 
= [1- \beta + \alpha \sqrt 3]/2. \label{eq:abtoxj}
\label{ew:xjtoab}
\enar
In terms of $\alpha$ and $\beta$,
Eq.~(\ref{eq:curve3}) for $\cC$ becomes,
\eq
3\alpha^2 - 2\sqrt 3 \mu \alpha\beta - 3\beta^2 + 2 \sqrt 3 \alpha
- 2 \beta +1 = 0, \quad \hbox{where} \ \
\mu = (k_1 - k_2)/(k_1 + k_2).
\label{eq:curveC}
\en
It is a hyperbola. Its intercept  with the $\beta$ axis
is given by the equation $3\beta^2 + 2\beta -1=0$,
which has one root,
$\beta = x_3  = 1/3$ while  $x_1 =x_2 = (1 - x_3)/2
=1/3$, representing the  point $E$ or its image $E^*$. 
The second  root, $\beta = -1$
should be excluded because $\beta = x_3$ cannot be negative. 
Thus the critical curve $\cC$ should be the upper branch
of the hyperbola (\ref{eq:curveC}) inside $\triangle_Q$, which is 
concave  upward in the $\alpha \beta$ plane and passes
through the centroid $E$ (or $E^*$ on the negative side).

The curve $\cC$   intersects the side $Q_3Q_1$ at 
point $\bx$ where $x_2=1/2= x_1 + x_3$,
and Eqs.~(\ref{eq:collinear})
and (\ref{eq:root3}) yield the  equations for $x_3$ and
$x_1 \in (0, 0.5)$,
\eq
x_3 = \frac{ \sqrt{k_1^2 + k_2^2 + k_1k_2}-k_2}{ 2 k_1}, \ \ \hbox{and}
\ \ \ x_1 = \frac{1}{2} - x_3 =
\frac{k_1 +k_2 -\sqrt{k_1^2 + k_2^2 + k_1k_2}}{2 k_1}.\ \
\label{eq:rootQ5}
\en
Likewise,  $\cC$ intersects  the side $Q_2Q_3$ at the critical point
$Q_4$ defined by Eqs.~(\ref{eq:collinear}) and (\ref{eq:root2}).

Thus $\cC$ is divided by $Q_4, E, Q_5$ and $E^*$ into four branches;
$Q_4 E$ and $EQ_5$ on the positive side of $\cP$ and
their images $Q_5E^*$ and $E^*Q_4$ on the negative side.
From Eq.~(\ref{eq:order}), we have $S_0 < 0$ when $A \not= 0$, 
hence $\dot p$ has the same sign
as $-\gamma C_0$, i.~e., the same sign as $\gamma (x_1-x_2)$. 
On the positive side of $\cP$, $\gamma = 1$, we have $\dot p < 0$ 
on the branch $EQ_5$ where $x_1 < x_2 < 1/2$ and $\dot p > 0$ on
the branch $Q_4E$, where $x_2 < x_1 < 1/2$. On the negative side,
$\gamma = -1$, the sign of $\dot p$ 
on the images of the branches is reversed.
Thus the points in  the branches $Q_4E$ and $E^*Q_5$, excluding
their end points, correspond to
expanding similar configurations $\triangle_p$ while those in
$EQ_5$ and $E^*Q_4$ correspond to contracting similar configurations.
It was observed before from  Eq.~(\ref{eq:soln-p}) that
the perimeter $p$ remains stationary at the end points
$Q_4, E, Q_5$ and $E^*$ and the corresponding vortex configuration
$\triangle_p$ or $\bR$ remains stationary also.

Now  we study the slopes and concavities 
of the trajectories, along constant $\bar I$ curves. From 
those  integral curves intersecting the critical curve $\cC$, 
we find the variation of the invariant
$\bar I$ along $\cC$ and then define the boundary of the strip
containing $\cC$. 
For this purpose, it suffices to study only  the positive side 
of $\cP$, where $\gamma = +1$.
We note that the slope of a point on the critical curve in the $\cP$
or $\alpha \beta$-plane, is defined by,
\eq
d\beta/d\alpha = (k_1x_2 - k_2 x_1)\sqrt 3/[2
(k_1+k_2)x_3 + k_2x_1 + k_1x_2].
\label{eq:slopeC}
\en
The slope has the same sign as the factor $(k_1x_2 - k_2x_1)$.
In the $\alpha \beta$ plane,
the curve $\cC$ is concave  upward with the minimum $\beta$ or $x_3$ at 
$ x_2 = k_2 x_1/k_1 \le x_1$. This is in agreement with
Eq.~(\ref{eq:curveC}), which  says that $\cC$
is the upper branch of a hyperbola passing through $E$.

For a trajectory in the parabolic case, Eq.~(\ref{eq:tt2}), 
its slope is,
\eq
m= d \beta/d\alpha=
\sqrt 3 x_3 (k_1 x_1-k_2 x_2)/[2(k_1+k_2)x_1x_2+(k_1x_1 + k_2x_2)x_3].
\label{eq:slopeI}
\en

The slope $m$ has the same sign as the factor $(k_1x_1 - k_2x_2)$, so
the trajectory  is concave downward with the maximum $\beta$ at
$ x_2 = k_1 x_1 / k_2 \ge x_1$ where $\alpha > 0$.

As $\beta$ increases from $0$ to $1/2$,
the invariant $\bar I$ decreases from $\infty$ monotonically 
along radial lines from $Q_3$, the origin, $\alpha=0, \beta = 0$.
Thus the  corresponding trajectories move upward as $\bar I$ decreases.

Since the critical curve $\cC$  is concave upward 
while the trajectories are concave downward, $\cC$  comes in contact
with the first trajectory when they are tangent to each other.
This happens  at point $E$ on $\cC$ with a trajectory with  $\bar I =1$
with  common slope
\eq
d \beta/d \alpha = (k_1 -k_2)/[(k_1 + k_2)\sqrt 3 ] \ge 0.
\label{eq:slope-E}
\en
That is, the trajectories with $\bar I > 1$
remain below $\cC$ and the one with $\bar I = 1$ is  tangent to $\cC$
at the  centroid $E$.
This statement holds for the images of the trajectories and the critical
curve on the negative side of $\cP$.
Hence {\it the trajectory with $\bar I = 1$ on both sides of $\cP$
is the lower boundary of a strip $\cal S$
containing $\cC$.
On the other hand, 
The trajectory  with $\bar I = 1$ encloses all the
trajectories  with $\bar I > 1$, which are the periodic orbits around
the center $Q_3$.}

As $\bar I$ decreases from $1$, the trajectory will intersect $\cC$
at two points on the left branch $Q_4 E$ and the right branch $EQ_5$.
respectively. To find the upper boundary trajectory of the strip $\cal S$,
we need to  find  whether it is the trajectory passing through 
$Q_4$ or  $Q_5$ by comparing the values of
their invariants $\bar I(Q_4)$ and $\bar I (Q_5)$. For this
comparison, 
we use $\bar I_j$ to denote $\bar I(Q_j)$.  Let the $x_3$ 
coordinates of $Q_4$ and $Q_5$ be denoted by $\beta_4$ and $\beta_5$,
then we have
\eq 
\bar I_4 =\bar I (Q_4) = (2 \beta_4)^{-(k_1 +1)}(1 -2 \beta_4)^{k_1}, 
\ \ \hbox{and} \ \  
\bar I_5 =\bar I (Q_5) = (2 \beta_5)^{-(k_1 + 1)}(1 -2\beta_5)^{k_1},
\label{eq:IQ4IQ5}
\en
with $k_2=1$ according to Eq.~(\ref{eq:order}).
From the $x_3$'s for $Q_4$ and $Q_5$, Eq.~(\ref{eq:rootQ5}), we have
\eq
2\beta_4 = (k_1 +1)/(q + k_1)\quad \hbox{and} \quad
2\beta_5 = (k_1 +1)/(q + 1),
\label{eq:beta4-5} 
\en
where $q= \sqrt {k_1^2 + k_1 +1}\ge \sqrt 3$ and the ratio,
\eq
f(k_1) = \bar I_4 /\bar I_5 =  
[(q-1)/(q+1)]^{k_1} [(q+k)/(q-k)] > 0.
\label{eq:I4I5}
\en
When $k_1 = k_2=1$, the integral curves are symmetric with respect to
the $\beta$-axis, the bisector of $\angle Q_3$. In particular, we have
$\beta_5 = \beta_4$, $f=1$ or $ \bar I_4 = \bar I_5$, and
the trajectory through $Q_4$ coincides with that through $Q_5$.

With $k_1 > k_2 =1$, we shall show the ratio, 
$f(k_1) =\bar I_4 / \bar I_5 \ge 1$ for $k_1\ge 1$.
With $f(1) =1$, we need only prove $df/dk_1 > 0$,
or rather $d \ln f /dk_1 > 0$. 
With $2 q d q/dk = 2k_1 +1$, we have
\eq
d \ln f/dk_1
=  -\ln [(1+\nu)/(1-\nu) + 2\nu + \nu (k_1^2 \nu^2)/(1 - \nu^2),
\label{eq:prf2}
\en
where $\nu = 1/q$. 
As $k_1$ increases from $1$ to $\infty$, 
$q$ increases from $\sqrt 3$ to $\infty$, 
$\nu$ decreases from $1/\sqrt 3$ to $0$ and  $(1-\nu^2)^{-1}$ decreases 
from $3/2$ to $1$. The curve $y = (1-x^2)^{-1}$ is concave upward and 
intersects the straight line or the chord, $y^+ = 1 + [x\nu /(1-\nu^2)]$,
at  points $(0, 1)$ and $(\nu, (1-\nu^2)^{-1})$. That is, the curve
lies below  the chord,
$y^+ > y$, for $x\in (0, \nu)$. 
Thus we have 
\eq
\frac{2-\nu^2}{2(1-\nu^2)}= \frac{1}{\nu} \int_0^{\nu} y^+ dx
> \frac{1}{\nu}\int_0^\nu \frac{d x}{1 -x^2} 
= \frac{1}{2\nu} \ln \frac{1+ \nu}{1 - \nu}
\label{eq:meanthm}
\en 
Equation (\ref{eq:prf2}) becomes 
\eqar
& & \quad d \ln f/dk_1  =
[(2 - \nu^2)\nu/(1-\nu^2) - \ln [(1+\nu)/(1-\nu)]
+ (k^2_1 - 1)\nu^3/(1-\nu^2) > 0,
\label{eq:prf3}
\\
& &\hbox{and hence} \quad
\bar I (Q_4) \ge \bar I (Q_5) \ \  \hbox{or} \ \ \bar I_4  \ge \bar I_5,
\quad \hbox{for} \ \ k_1 \ge 1
\label{eq:prf5}
\enar
The equality sign holds only when $ k_1 = 1 =k_2$. 
This completes the proof,
which was omitted in the 1988 paper \cite{tavtin}.
Note that  the trajectory through $Q_4$ lies in the strip
$\cal S$, above the lower trajectory  through $E$ and below
the upper one   through $Q_5$ as shown in Fig.~\ref{fig:parabA}. 

For the trajectories with $\bar I \in (\bar I_6 , \bar I_5)$, 
they will remain above 
$\cC$ and cross over to the opposite side of $\cP$ to form periodic orbits.
For $\bar I = \bar I_6$, the trajectory  becomes the separatrix
through the critical point $Q_6$ on the side $Q_1 Q_2$, 
as shown in Fig.~\ref{fig:parabA}. For $\bar I \in (0, \bar I_6)$
the trajectories will either be periodic orbits around the center $Q_1$
or the center $Q_2$. This completes the description of the trajectories,
separatices and the critical curve
in a parabolic case shown in
Fig.~\ref{fig:parabA}, which was presented in \cite{tavtin}.

From Eqs.~(\ref{eq:prf5}) and (\ref{eq:slope-E}), we see that
{\it on the positive side of $\cP$,
the critical curve $\cC$ is contained in a strip $\cS$, bounded
above  by the trajectory  through point $Q_5$ on edge $Q_3 Q_1$
with $\bar I =\bar I_5$ and below by the trajectory through 
the centroids, $E$ and $E*$,
with $\bar I = 1$. Thus the strip $\cS$
contains besides $\cC$ all trajectories with $\bar I \in [\bar I_5, 1]$.
Likewise, we define the image  of the strip on 
the negative side of $\cP$.}

In the next subsection we will find out whether 
a trajectory through a point near $\cC$ will be attracted to 
a  point $\bq$ on $\cC$
or depart from $\bq$  with increasing time.

\subsection{Trajectories in the strip $\cS$ containing $\cC$}
\label{ssec-direc}

Note that a point on the critical curve $\cC$ is a stationary point.
A trajectory will either stay above $\cC$ or below $\cC$ 
but {\it will not cross over $\cC$}. Thus $\cC$ partitions
the strip $\cS$ on the positive side of $\cP$ into
two strips $\cS^+$ and $\cS^-$ above and below $\cC$.
To find out whether a trajectory through a point in the neighborhood 
of $\cC$ will depart from or be attracted to $\cC$, we need to know
the direction
along the trajectory as $t$ increases. 
For this, we use Eqs.~(\ref{eq:dotxj}).
We see that $\dot x_1$ and $\dot x_2$ changes signs only when crossing
over an edge where $\gamma$ changes sign or over $\cC$ where $\cal Y$
does. The factor $F_1 > 0$  and $F_2 < 0$ so that they contribute
to a factor in $\dot \alpha < 0$ because it has the same sign as that
$F_2 - F_1$. That is the $\dot \alpha$ will have the same sign
as $\dot x_2$ and change its sign in crossing over an edge or over $\cC$.
For $\dot x_3$ or $\dot \beta$,
we see there is an addition sign change because 
the factor $F_3$ can cross over zero. 
This sign change in $\dot \beta$ defines the concavity of
the trajectory, concaving downward, and the maximum of $\beta$.
Thus the directivity of the trajectories are defined by
Eq.~(\ref{eq:dotxj}).
For $j=3$ we  obtain
\eq
\dot \beta = \dot x_3 = \gamma {\cal H} {\cal Y} x_3 [\frac{k_2 x_2
-k_1 x_1}{k_1k_2}].
\label{eq:direction}
\en
and $\beta$ or $x_3$ reaches its maximum when $k_2x_2= k_1x_1$ as expected from
Eq.~(\ref{eq:slopeI}) for the slope $m=0$. From the concavities of the
trajectories and the critical curve, we showed in
the preceding subsection that a trajectory
with $\bar I \in [\bar I_4, 1)$ will intersect $\cC$ twice in $\cS^+$ on
the positive side of $\cP$ and cross over the two edges in $\cS_-$
at the expanding branch, $Q_4E$ to the left of the centroid $E$, 
where $x_1 > x_2$ with slope $m > 0$ and at the contracting  branch  
$EQ_5$ to the right of $E$ where $x_2 > x_1$ with $m < 0$. 
A trajectory with $\bar I \in [\bar I_5, \bar I_4)$  will
intersect only the right branch $EQ_5$ of $\cC$  and
cross over the left edge in $\cS^+$ and the right edge in $\cS^-$.
From Eq.~(\ref{eq:dotxj}), we see that 
on the positive side of $\cP$, $\gamma = 1$, and $\dot \beta$
has the same  sign as ${\cal Y} [k_1x_1 - k_2x_2]$ or the opposite
sign of $m {\cal Y}$ since ${\cal H} >0$.

Recalling  that  ${\cal Y} >0$ when the point is above
$\cC$, and ${\cal Y} < 0$ when below, we have the direction of
a trajectory near $\cC$ on the positive side of $\cP$, 
\eqar
&\hbox{Near the branch}\ EQ_5, \quad m > 0, \qquad \qquad&  \dot \beta > 0 
\ \ \hbox{above}\ \cC ,  \\
&        & \dot \beta < 0, \ \ \hbox{below}\ \cC, \\
&\hbox{Near the branch}\ Q_4E, \quad m < 0, \qquad \qquad&  \dot \beta < 0 
\ \ \hbox{above}\ \cC ,  \\
&        & \dot \beta > 0, \ \ \hbox{below}\ \cC.
\enar
On the negative side of $\cP$, $\gamma = -1$, the inequality signs will
be reversed. The directions of the trajectories shown in
Fig.~\ref{fig:parabA} are in agreement with the statements above.
These statements say that

\noindent{\it In the neighborhood of a contracting (expanding)
branch, $EQ_5$ or $Q_4 E^*$
($Q_4E$ or $E^*Q_5$), a trajectory will depart from (be attracted to)
$\cC$.}

It is equivalent to say that

\noindent{\it When the initial vortex configuration similar to
$\triangle_1$ or point $\bq$ 
on a contracting branch of
the critical curve $\cC$ is disturbed off $\cC$, the trajectory 
will depart from $\cC$ and finally be attracted to an expanding branch
of $\cC$.}

In particular,  
\noindent{\it Contracting similar vortex configuration, $\triangle_p$
leading to coelescence of the three vortices is a theoretical
solution but unlikely to happen.}

Knowing the directivities of the trajectories, we shall describe 
in the following subsection how the 
trajectories in the strip $\cS$ 
depart from a contracting or unstable branch of
$\cC$ ending on an expanding or stable branch
and classify the trajectories  into three types.

\subsection{Departure from a contracting branch of the critical curve}
\label{ssec-bifurcation}

For the region $\cS^-$, the centroid $E$ divides it into two with
$x_2 > x_1$ and $x_2 <  x_1$ respectively. Likewise the image of
$\cS_-$ is divided by $E^*$ into two.
We use $\cS^-_R$ to denote the  subregion of $\cS^-$ where $x_2 > x_1$,  
bounded 
above by the contracting branch $EQ_5$ and below by the trajactory
with $\bar I = 1$ on the positive side of $\cP$
and is connected across the edge where $x_2 = 0.5$ to
its image on the negative side of $\cP$ in which the upper
boundary is the expanding branch $Q_5E^*$.
A trajectory in $\cS^-_R$ with $\bar I \in (\bar I_5, 1)$
will depart from a point below the
contracting branch $EQ_5$ cross over to edge and be attracted to
the expanding edge $E^*Q_5$.
Likewise, we use $\cS^-_L$ to denote the corresponding region where
$x_2 < x_1$. The trajectories 
$\cS^-_L$ have $\bar I \in  (\bar I_4, 1)$ would depart
from the contracting branch $E^*Q_4$ and be attracted to the expanding
branch $Q_4E$.
The trajactories in $\cC^-$
are classified as trajectories of type I.

For the region $\cS^+$ above $\cC$,
Eq.~(\ref{eq:prf5}) says that  
there is a separatrix from $Q_4$
which intersects $\cC$ at $S_4$ and partitions the branch $EQ_5$
into $ES_4$ and $S_4 Q_5$ and also the strip $\cS^+$
into two strips $\cS^+_4$ and $\cS^+_5$ 
with $\bar I \in (\bar I_4, 1)$ and $\in (\bar I_4, \bar I_5)$,
respectively
(see Fig.~\ref{fig:parabA}).

For a $\bar I \in (\bar I_4, 1)$, the trajectory in $\cS^+_4$
on the 
positive side of $\cP$ intersects $\cC$ at 
two points on  the left branch $Q_4 E$ and the right branch
$ES_4$, respectively, i.~e., it
remains  on the positive side.
The same holds for its image on
the negative side. The trajectory in $\cS^+_4$ 
and its image are of type II (see Fig.~\ref{fig:parabA}).

\begin{figure}[hbt]
\centerline{
\includegraphics[angle=0,width=5.0in, draft=false]{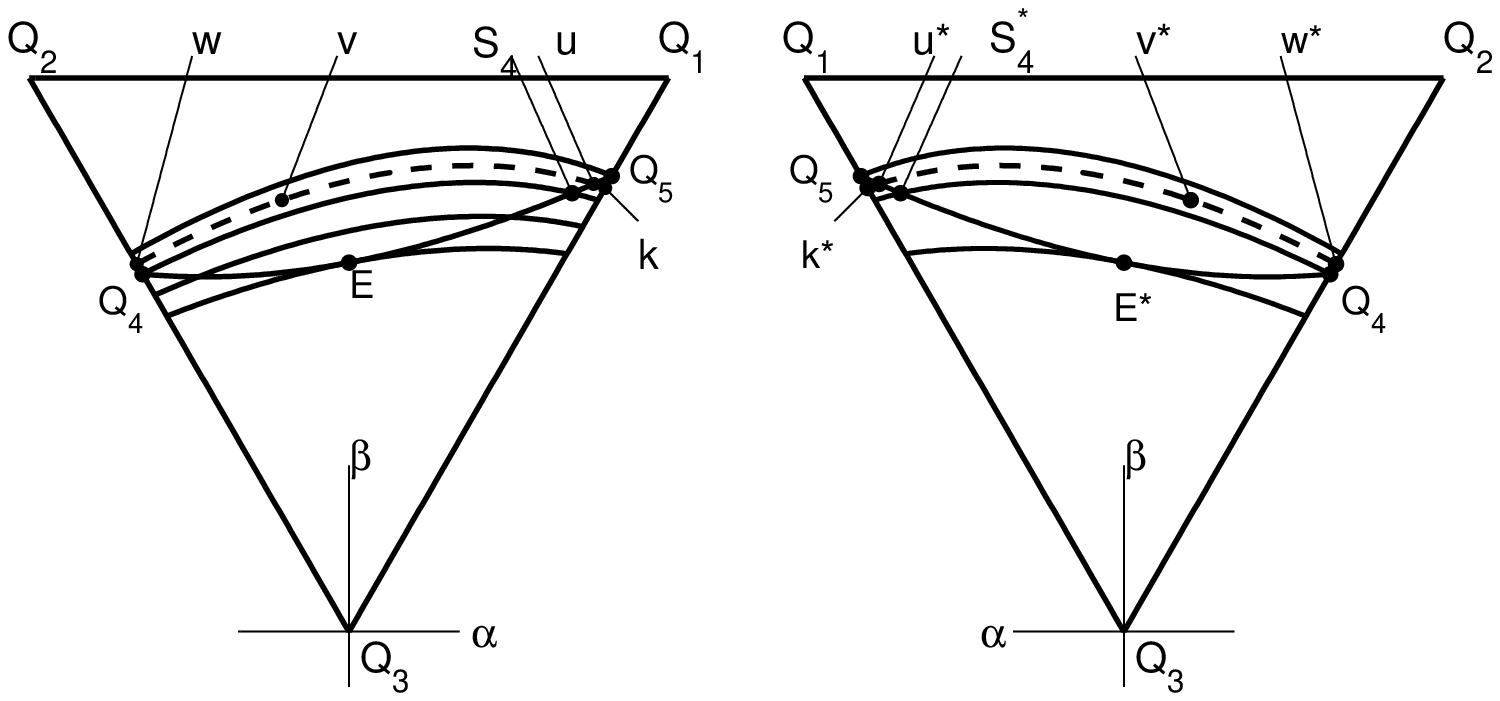} }
\caption{Departure from a  point, $u$ or $u^+$,   on
an unstable branch $EQ_5$ of the critical curve 
along the trajectory $u - v - w$ in the plane $\cP$, with $\gamma =+1$
and along its image $w- v^* - u^* - u_c$ in $\cP$ with $\gamma = -1$.
The last point $u_c$ denotes the point  where the trajectory 
intersects the stable branch $E^*Q_5$ of $\cC$.} 
\label{fig:front2}
\end{figure}

For a $\bar I \in (\bar I_5, \bar I_4)$, the trajectory in
$\cS^+_5$ in the positive part of $\cP$ will be above $\cC$,
cross over the edge $Q_2Q_3$ and continue along the image of 
the trajectory in the reverse direction. The trajectories in
$\cS^+_5$ are classified as type III (see Fig.~\ref{fig:front2}).

The above theoretical predictions and classifications
of the trajectories in the strip $\cS$ will be illustrated or
confirmed by the numerical examples in the next section.

\section{Numerical Simulations}
\label{sec-numerical}

In this section, we present numerical examples to simulate
or substantiate the conclusions arrived 
at in  Subsec.~\ref{ssec-bifurcation}
via the directivity of the trajectories
intersecting the critical curve $\cC$. That is, a configuration
or a point $\bx$ in the neighborhold of a contracting branch of 
of $\cC$ will depart from $\cC$ along a trajectory with constant
integral invariant $\bar I$ and end on
an expanding branch of $\cC$.
Instead of 
integrating the govening equations~(\ref{eq:dotRj})
for the vortex configurations from the initial configuration,
the numerical examples are based on the integration
of the primary
equations of motion of the three
vortices, Eqs.~(\ref{eq:motion}) in the $z$-plane with initial data 
$z_j(0)$'s corresponding to the initial
configuration $\triangle_1$'s in the neighborhood of $\cC$. 

First we prescribe the strengths of
the three vortices belonging to the parabolic case: with
$k_2= 1$ as the scale for the strength, we choose $k_1 = 2$
and get $k_3 = k_1 k_2/(k_1 + k_2)= -2/3$. 
Without losing generality, 
we choose the initial perimeter of $\triangle_p$ to be 
the length scale,
i.~e., $p(0) = 1$ and hence $R_j(0) = x_j(0)$. 

To carry out the integration of
the primary equations~(\ref{eq:motion}), we shall relate the
initial positions of the vortices\footnote{Note that 
$x_j, \ j =1, 2, 3$ denote the trilinear coordinates and  should not 
be confused with $\Re z_j$.}
$z_j(0)$'s to  the initial configuration, 
$R_j(0) = x_j(0)$  by
(i) aligning $z_1$ and $z_2$ parallel to the imaginary axis, i.~e., 
setting $\Re z_1 = \Re z_2 $,
(ii)  requiring $\Im z_1 > \Im z_2$ and $\Re z_3 > \Re  z_1$, so that
the configuration $\triangle_1$ is oriented counterclockwise,
$\gamma = +1$, and
(iii) putting the stationary center of vorticity  at the origin, $z=0$.
The relationships for the initial data with $R_j(0) = x_j(0)$ are:
\eqar
& & |z_3 - z_1|^2  = R_2^2 \qquad 
 |z_2 - z_3|^2 = R^2_1 \label{eq:R2R1} \\
& & z_1  - z_2  = i R_3 \qand
  k_1z_1+ k_2 z_2 + k_3 z_3 = 0.
\label{eq:kjzj}
\enar
From Eq.~(\ref{eq:R2R1}) and the imaginary parts of the equations in
(\ref{eq:kjzj}), we get  $\Im z_j(0)$,
\eqar
& & \Im z_2  = [-k_3 (R_1^2 - R_2^2) 
+R_3^2 (2 k_1 + k_3)]/[2R_3(k_1+k_2+k_3)], 
\label{eq:Imz2} \\
& & \Im z_1 = \Im z_2 + x_3 \qand \Im z_3
= - [k_1 \Im z_1 + k_2 \Im z_2]/k_3 
\label{eq:Imz13}
\enar
By equating the real parts of the equations in (\ref{eq:kjzj}), we get
\eq
\Re z_1 = \Re z_2, \ \  \Re z_3 - \Re z_1 = 
\sqrt{R_2^2 - (\Im z_3-\Im z_1)^2} \quad     
\hbox{and} \ \ (k_1 + k_2) \Re z_1 + k_3 \Re z_3 = 0.
\label{eq:xjsol}
\en
These three equations in turn define $\Re z_j(0)$.
With those initial data, we carry out 
numerical integration of the primary equations~(\ref{eq:motion}),
obtain the motions of the vortices, $z_j (t)$'s,
and  recover  the vortex configuration,
the sides $R_j(t)$ and the  perimeter $p(t)$ of $\triangle_p (t)$,
and then the three sides, $x_j(t)= R_j/p$ of $\triangle_1(t)$. 

Ten  numerical examples were carried out and the results
concurred with the theoretical conclusions that
a configuration slightly  perturbed from a contracting branch of the 
critical curve $\cC$
will depart from $\cC$ and  end  
on  an expanding branch of
$\cC$ belonging to one of the three types described in
Subsec.~\ref{ssec-bifurcation}. 
Here we shall present only four typical  numerical examples.
Table 1 gives the
four sets  of  initial configurations, 
$R_j(0)= x_j(0)$ with $p(0) = 1$ and the corresponding 
values of $Y(0)$ 
and the invariant $\bar I(0)$ defined by Eqs~(\ref{eq:curve3})
and (\ref{eq:tt1}), respectively.

\begin{table}
\begin{center}
\begin{tabular}{|c|c|c|c|c|c|}\hline
Case &  $R_1$   & $  R_2$  &  $R_3$   &    $Y(0)$  &$\bar{I}(0)$\\ 
\hline
$r_-$   & 0.18195  &  0.44396 &  0.37409 & -0.00498 &  0.68503 \\
$r_+$   & 0.19108  &  0.43424 &  0.37468 &  0.00501 &  0.68500 \\
$u_-$   & 0.10442  &  0.49225 &  0.40333 & -0.00500 &  0.38563 \\
$u_+$   & 0.10839  &  0.48643 &  0.40518 &  0.00502 &  0.38555 \\
\hline
\end{tabular}
\caption{Specification of numerical examples.  
In all cases, $k_1=2$, $k_2=1$, and $k_3 = -2/3$.}
\label{table:1}
\end{center}
\end{table}

Since the critical curve $\cC$
is defined by $Y(x_j) = 0$, see Eq.~(\ref{eq:calY}), 
the deviation of point $x_j$ from
$\cC$ is measured by $Y(x_j)$ with $Y$ positive (negative)
for $x_j$ lying  above (below) the critical curve. 

The contracting branch $EQ_5$ of $\cC$ on the positive side of the
plane $\cP$ is composed of segments $ES_4$ and $S_4 Q_5$ with invariant
$\bar I$ decreasing respectively from $\bar I(E) = 1$ to $\bar I_4$
and from $\bar I_4$  to $\bar I_5$.
See Figs. \ref{fig:parabA} and \ref{fig:front2}.

The integral invariants $\bar I (0) $'s  in the table
define the locations
of the critical points on the contracting branch of $\cC$ 
and each shall remain constant for $t >0$
along the trajectory towards an expanding branch of $\cC$.

The  first and second rows of data in Table 1
list the initial data corresponding to point $r^{\mp}$ perturbed 
below and above, respectively. from point $r$ in the segment
$E S_4$ of $\cC$.
Likewise, the third  and fourth  rows 
list  the  initial data corresponding to points $u^{\mp}$ 
perturbed from a point $u$ 
in the segment $S_4Q_5$ of $\cC$.\footnote{The $x_j$'s of point $r$ 
on $\cC$, where $Y=0$ can be defined from the linear 
interpolation of $x_j$'s  of $r^{\pm}$ with respect to $Y^{\pm}(0)$.
Likewise we can locate the  $x_j$is  of point $u$ on $\cC$
from those of $u^{\pm}$.}

For each case, the integral invariant $\bar I(t)$ differs from 
its initial value
$\bar I(0)$ by less than $0.1 \%$, that is, the 
numerical solution $z_j(t)$ does yield
a planar trajectory of constant $\bar I(0)$.

\begin{figure}[htb]
\centerline{
\includegraphics[angle=0,width=64mm,draft=false]
{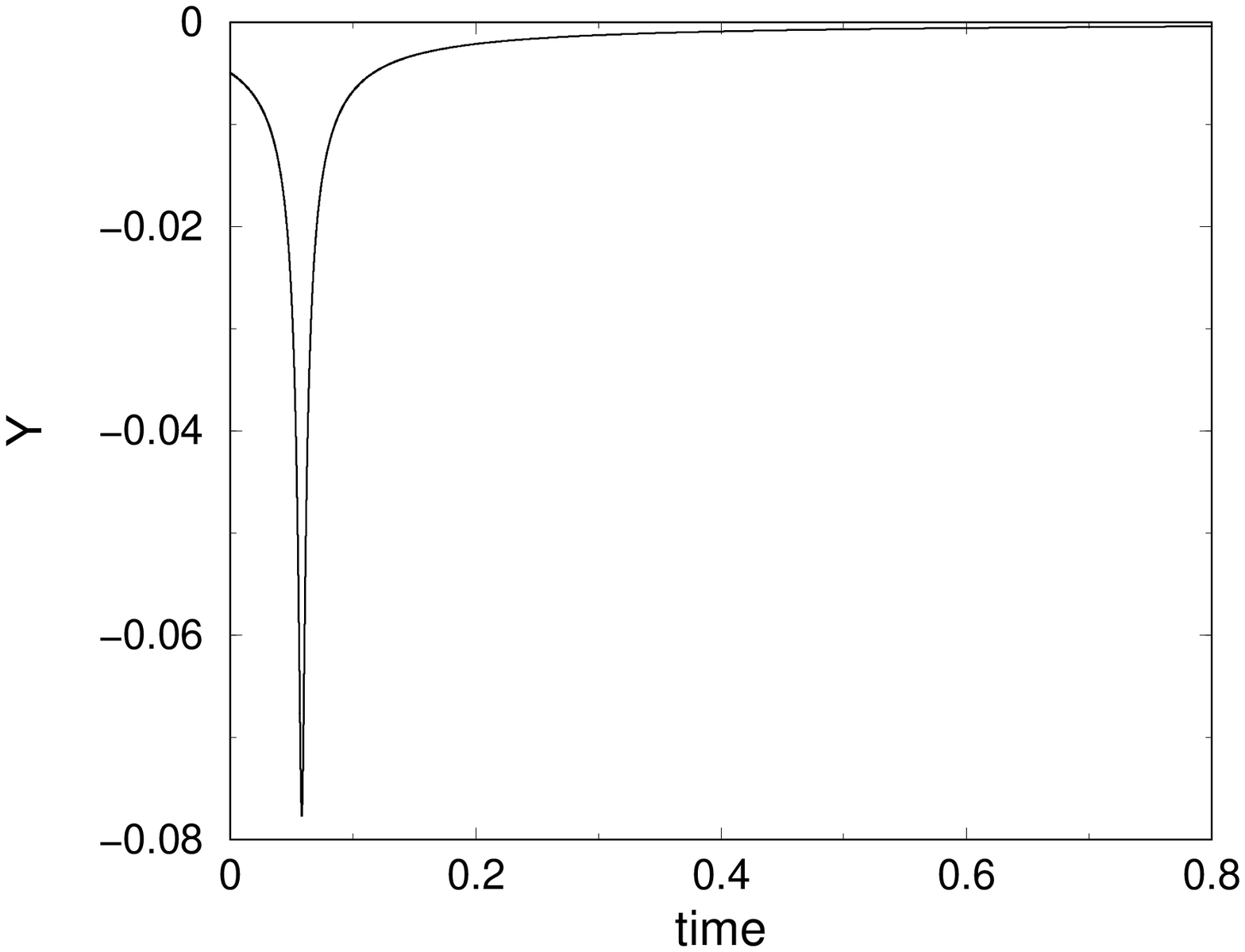}~a) \hspace{1.0mm}
\includegraphics[angle=0,width=64mm,draft=false]
{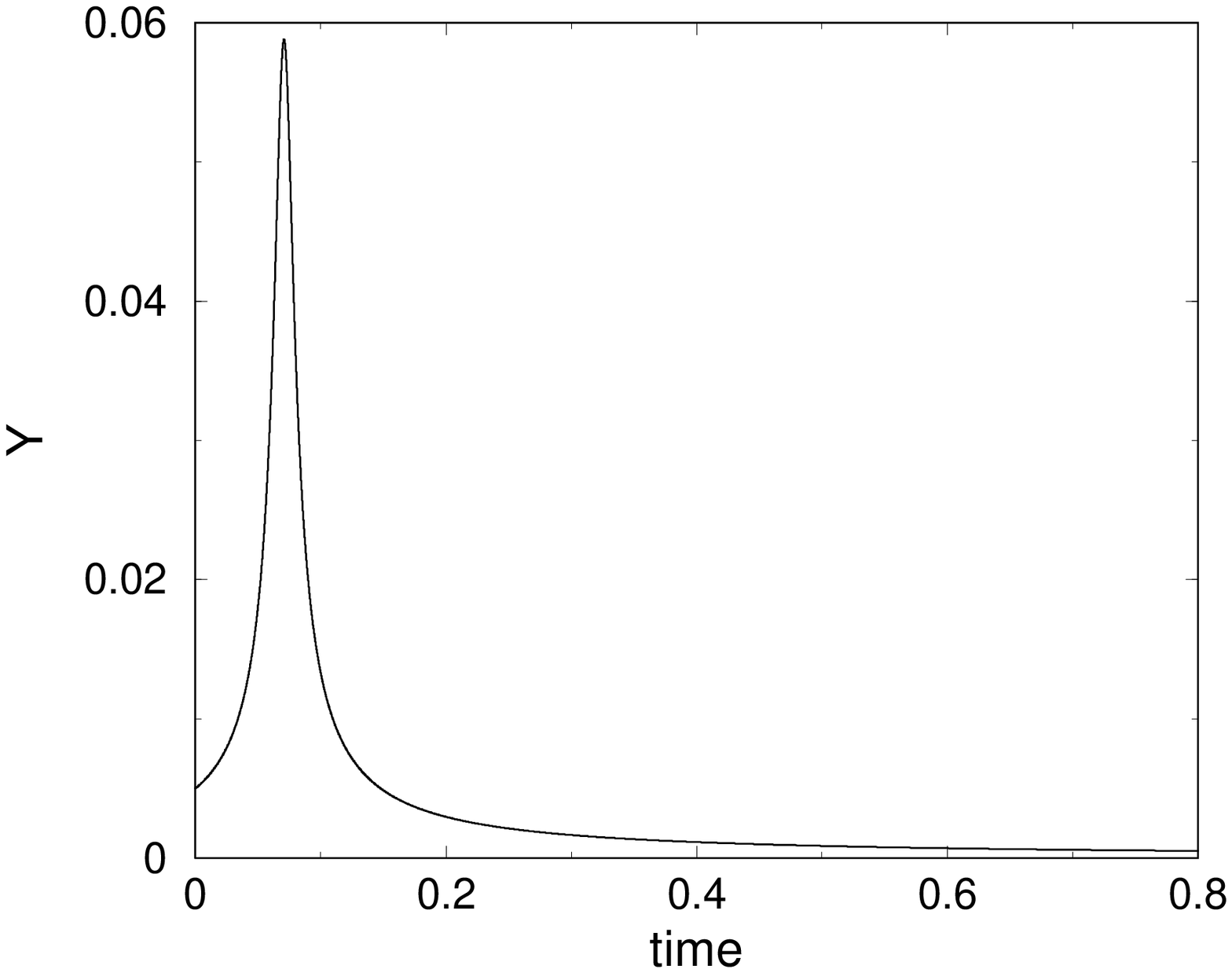}~b) }
\centerline{
\includegraphics[angle=0,width=64mm,draft=false]
{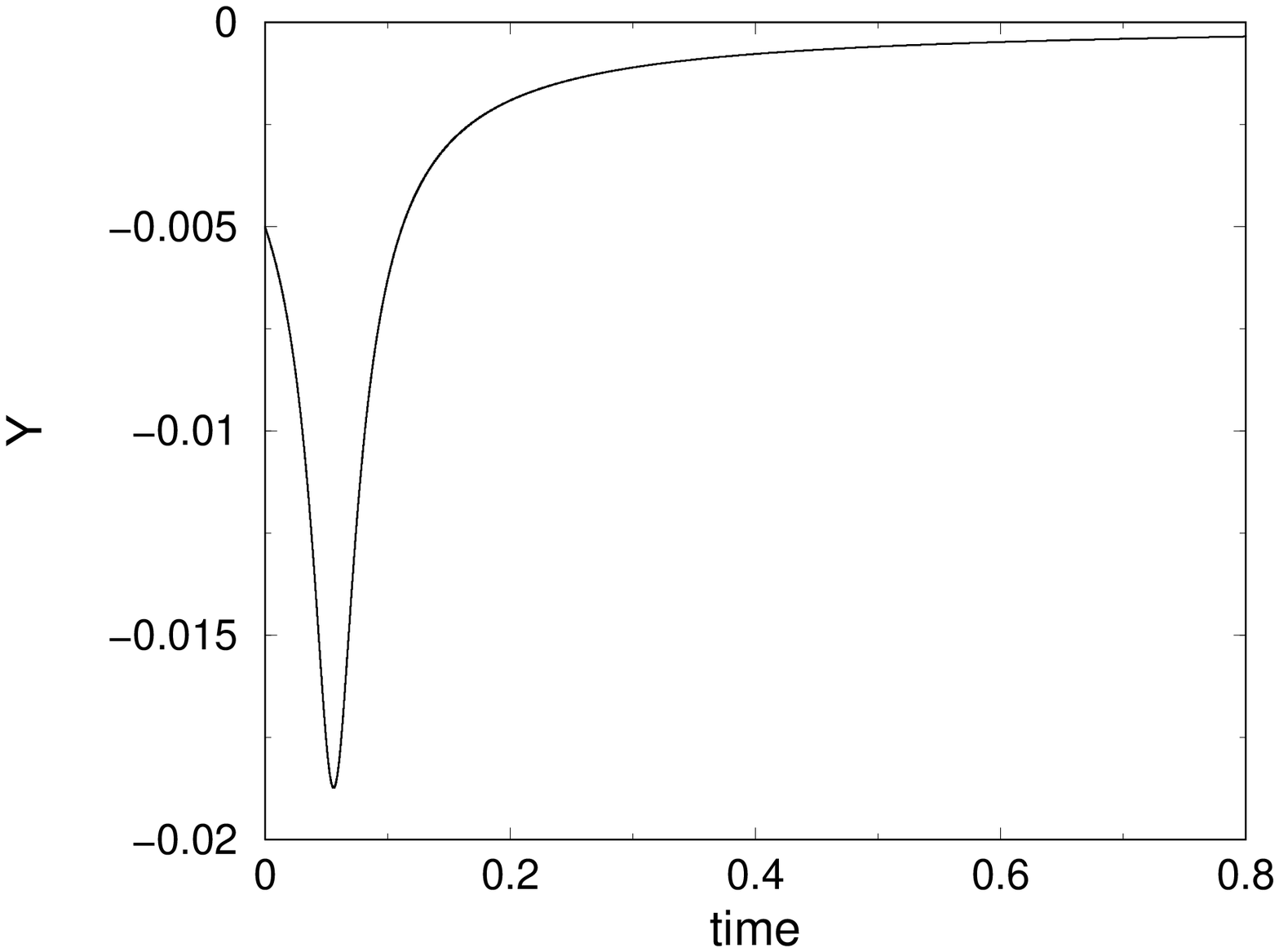}~c) \hspace{1.0mm}
\includegraphics[angle=0,width=64mm,draft=false]
{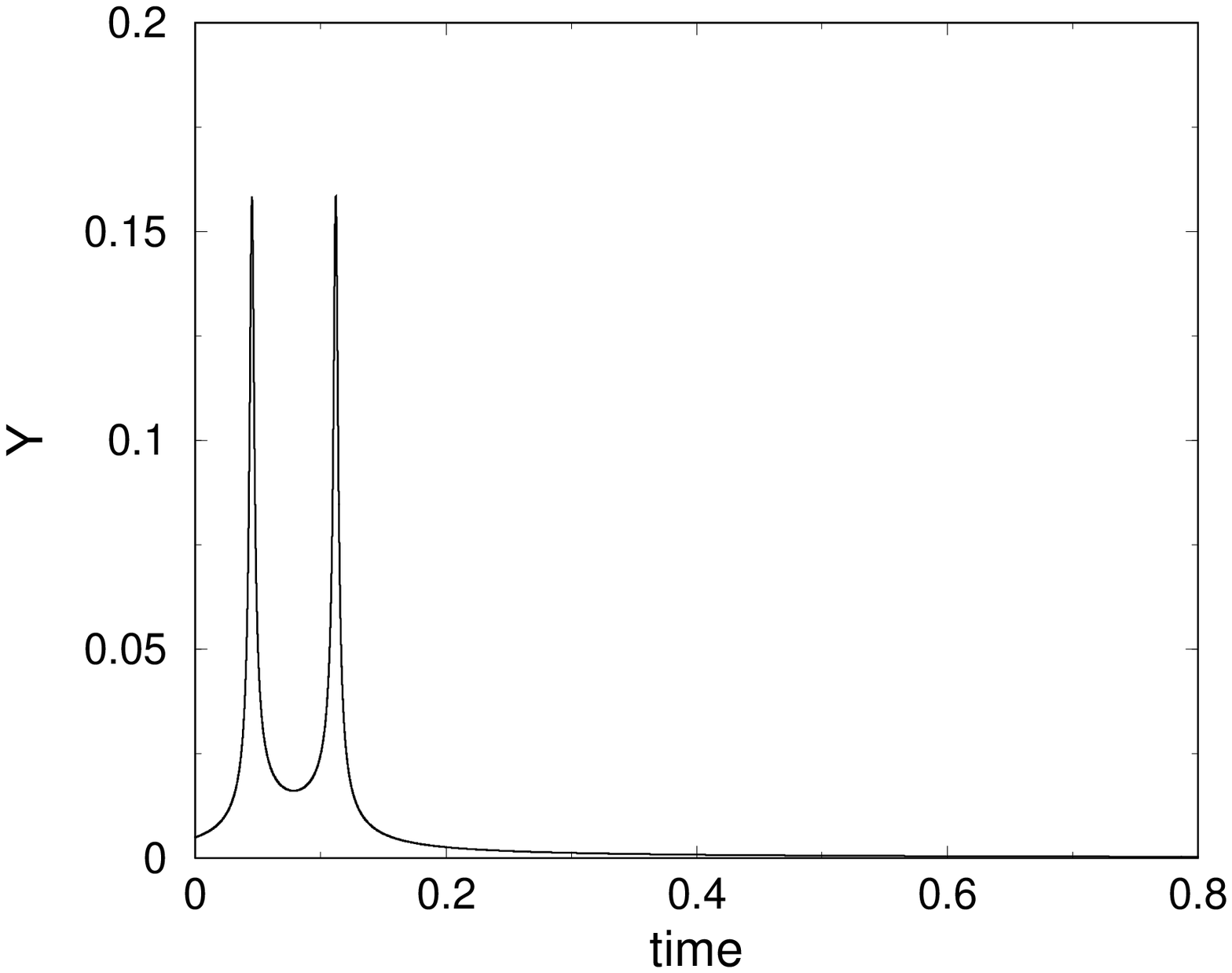}~d)}
\caption{Variation of $Y({\bf x})$ along a trajectory, ${\bf x} (t)$,
departing from
an unstable branch of the critical curve to a stable branch.
Point of departure: $r^-$, $r^+$, $u^-$, and $u^+$ in a),  b),  c) 
and  d).}
\label{fig:twoby2-Y}
\end{figure}

\begin{figure}[htb]
\centerline{
\includegraphics[angle=0,width=60mm,draft=false]
{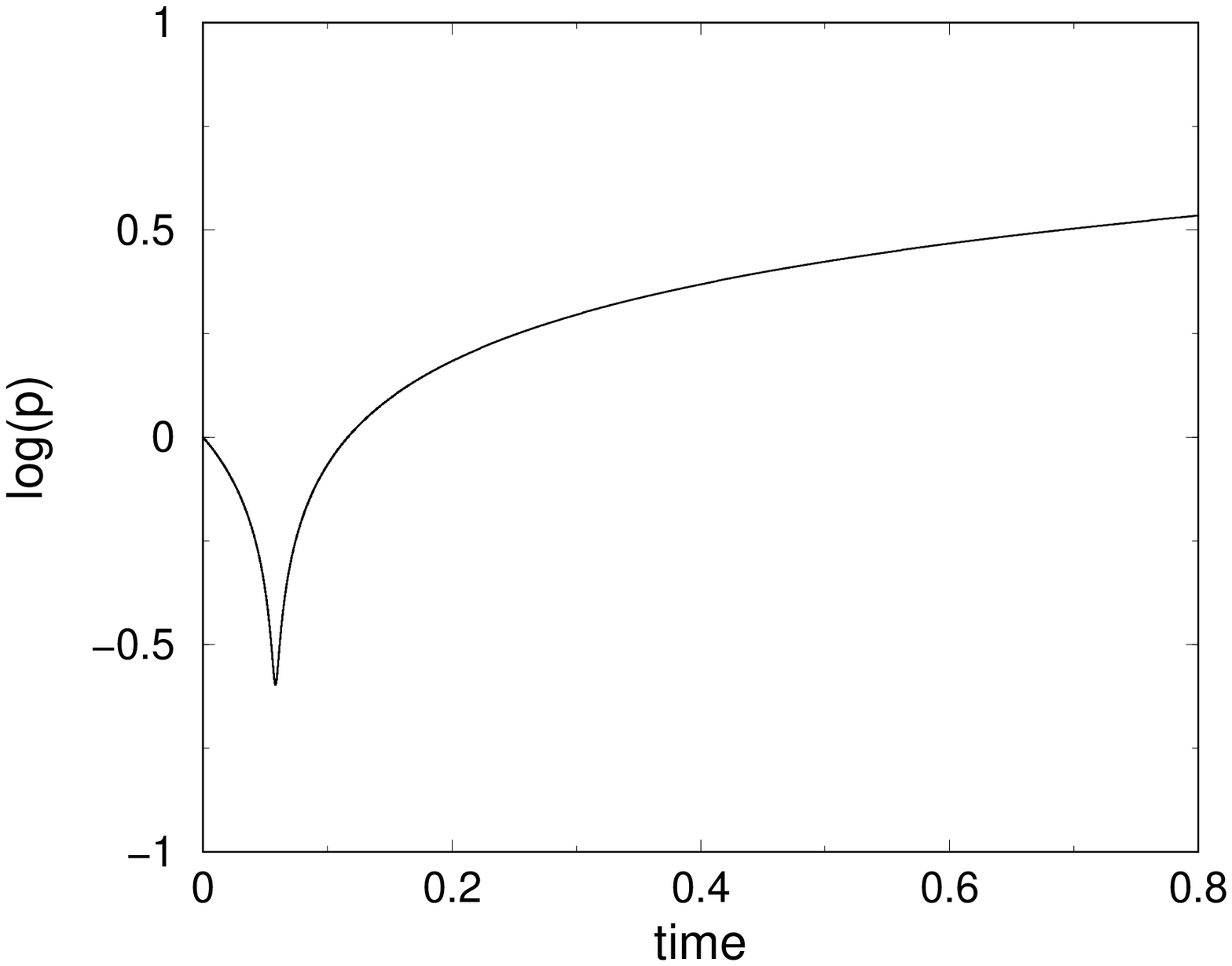}~a) \hspace{1.0mm}
\includegraphics[angle=0,width=60mm,draft=false]
{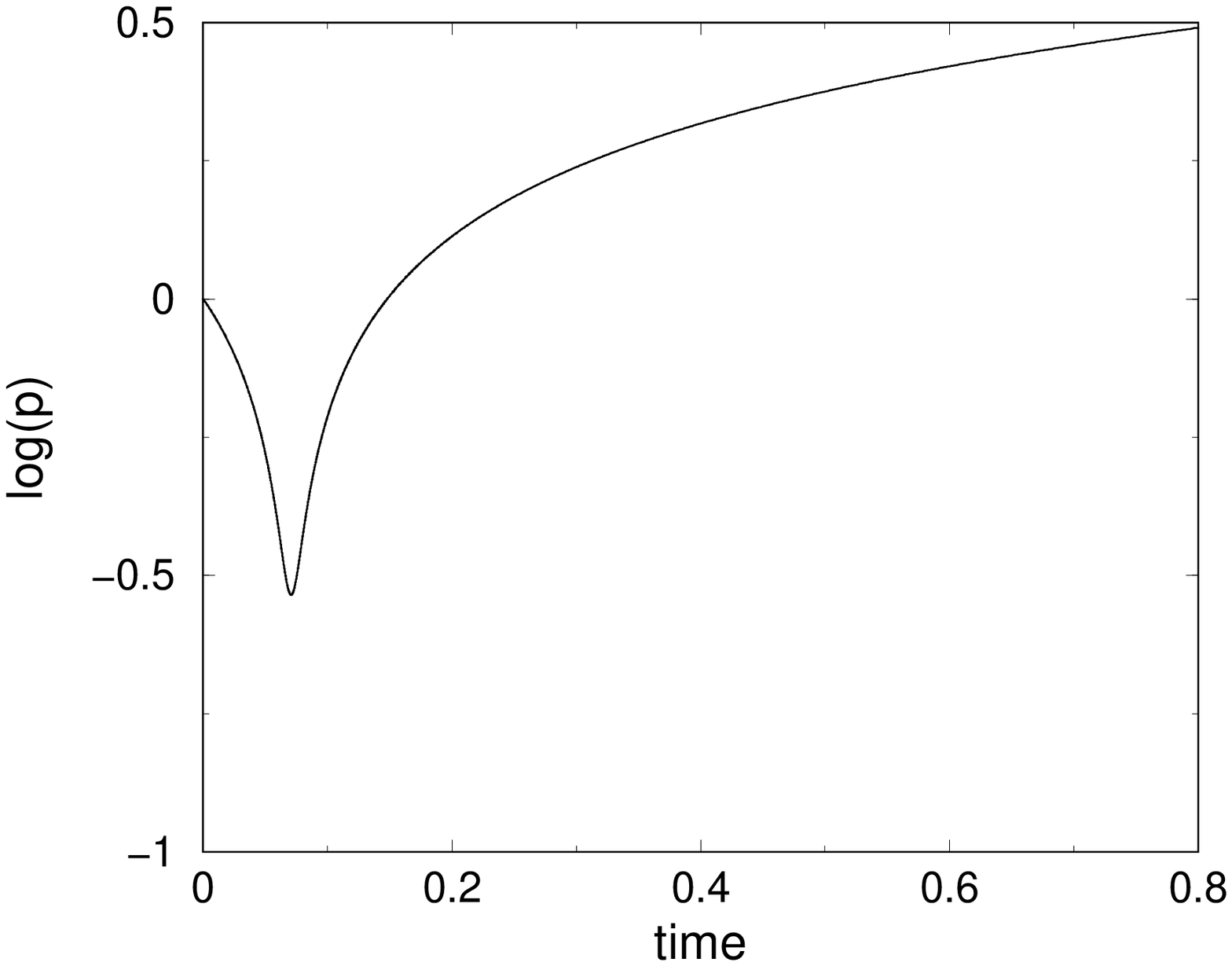}~b) }
\centerline{
\includegraphics[angle=0,width=60mm,draft=false]
{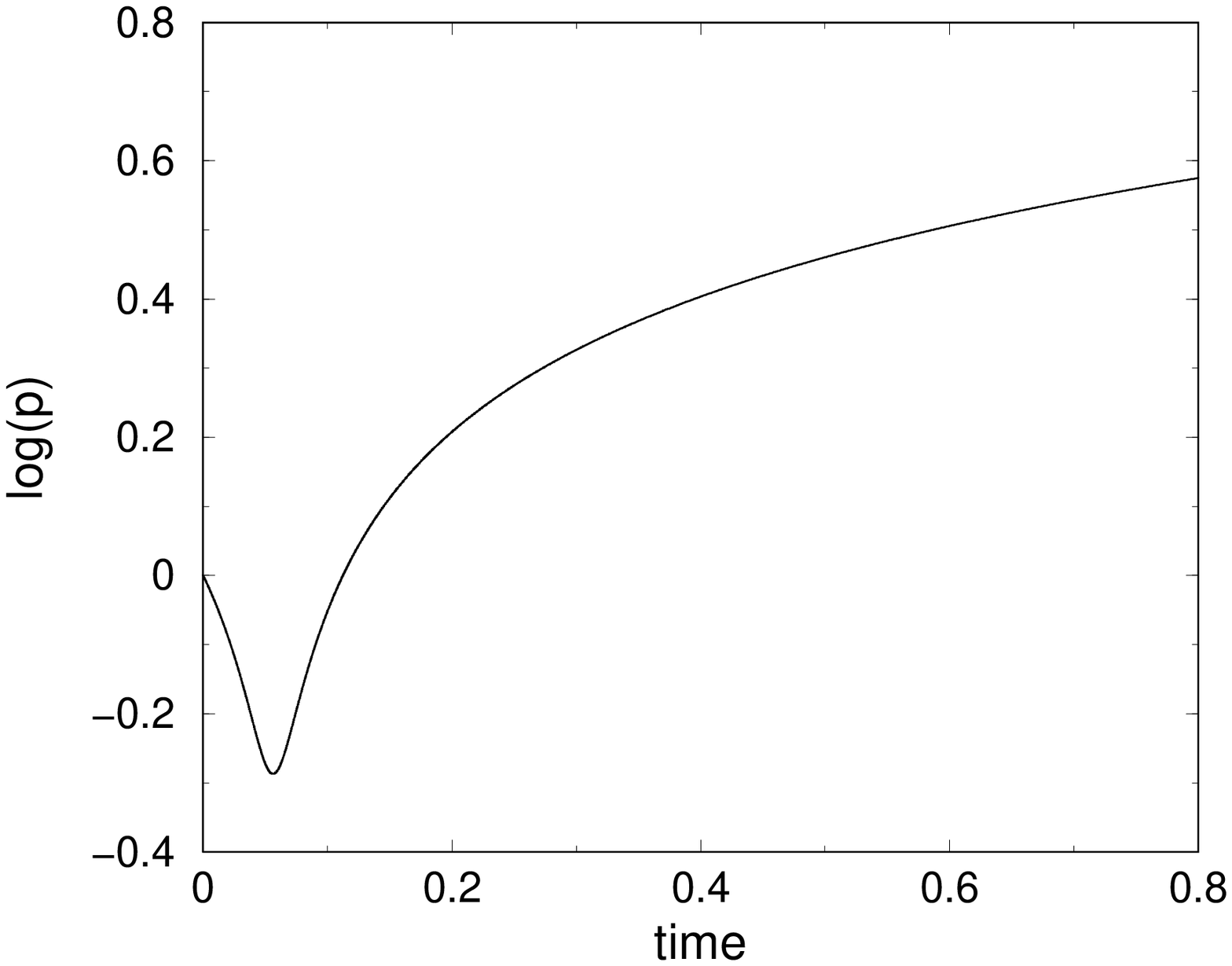}~c) \hspace{1.0mm}
\includegraphics[angle=0,width=60mm,draft=false]
{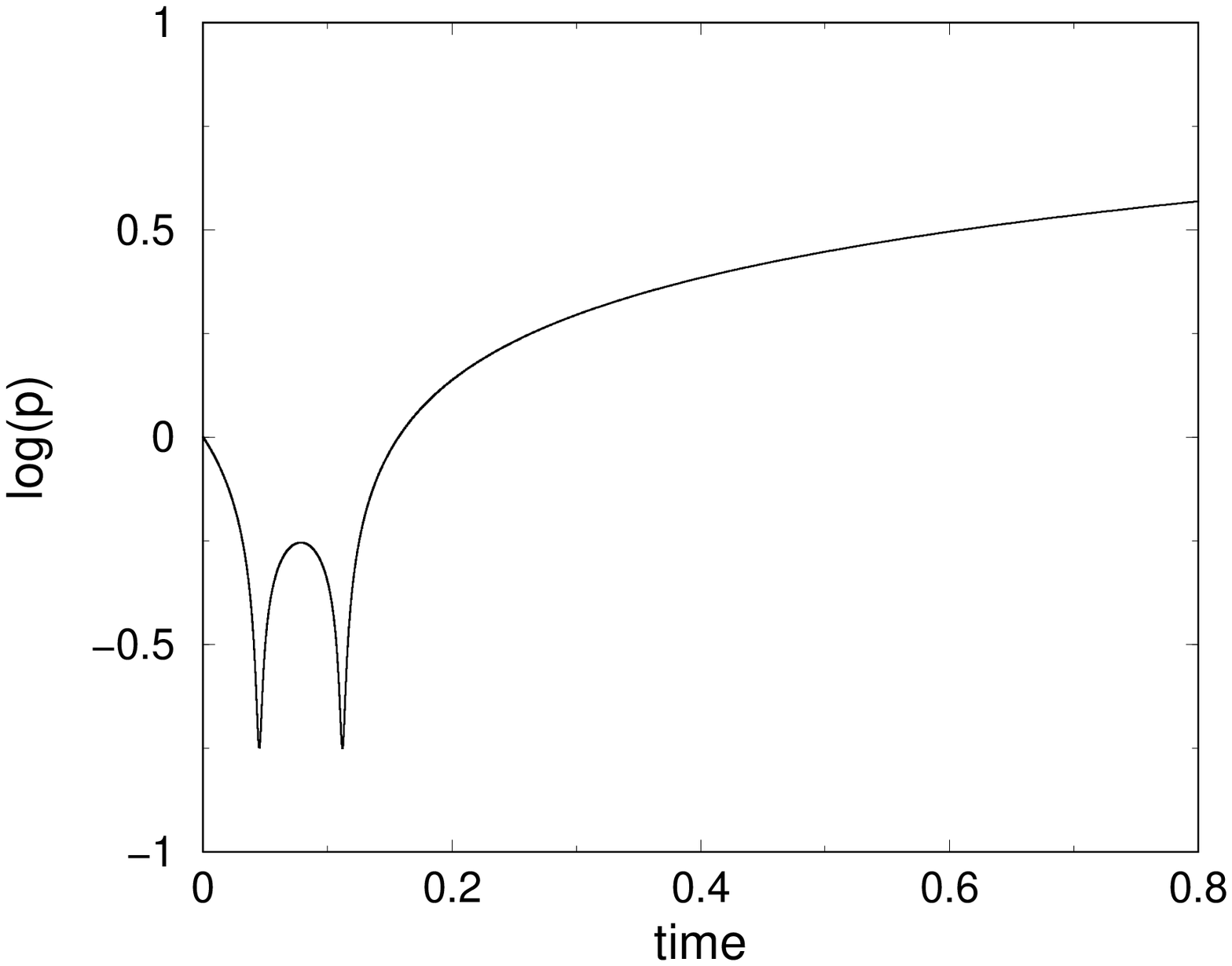}~d)}
\caption{Variation of $\log p({\bf x})$ along a trajectory, 
${\bf x} (t)$, departing from
an unstable branch of the critical curve to a stable branch.
Point of departure: $r^-$, $r^+$, $u^-$, 
and $u^+$ in a),  b),  c) and  d).}
\label{fig:twoby2-p}
\end{figure}

\begin{figure}[htb]
\centerline{
\includegraphics[angle=0,width=4.5in,draft=false]
{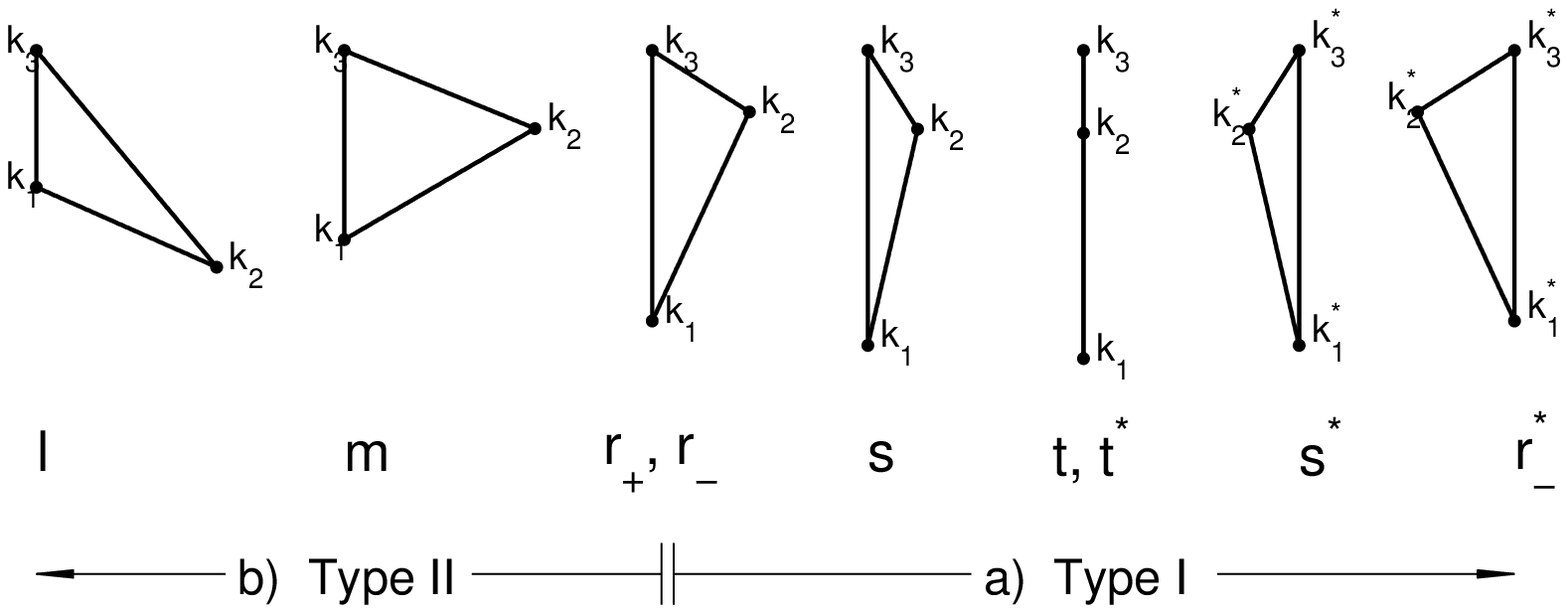}} 
\vspace{0.25in}
\centerline{
\includegraphics[angle=0,width=4.5in,draft=false]
{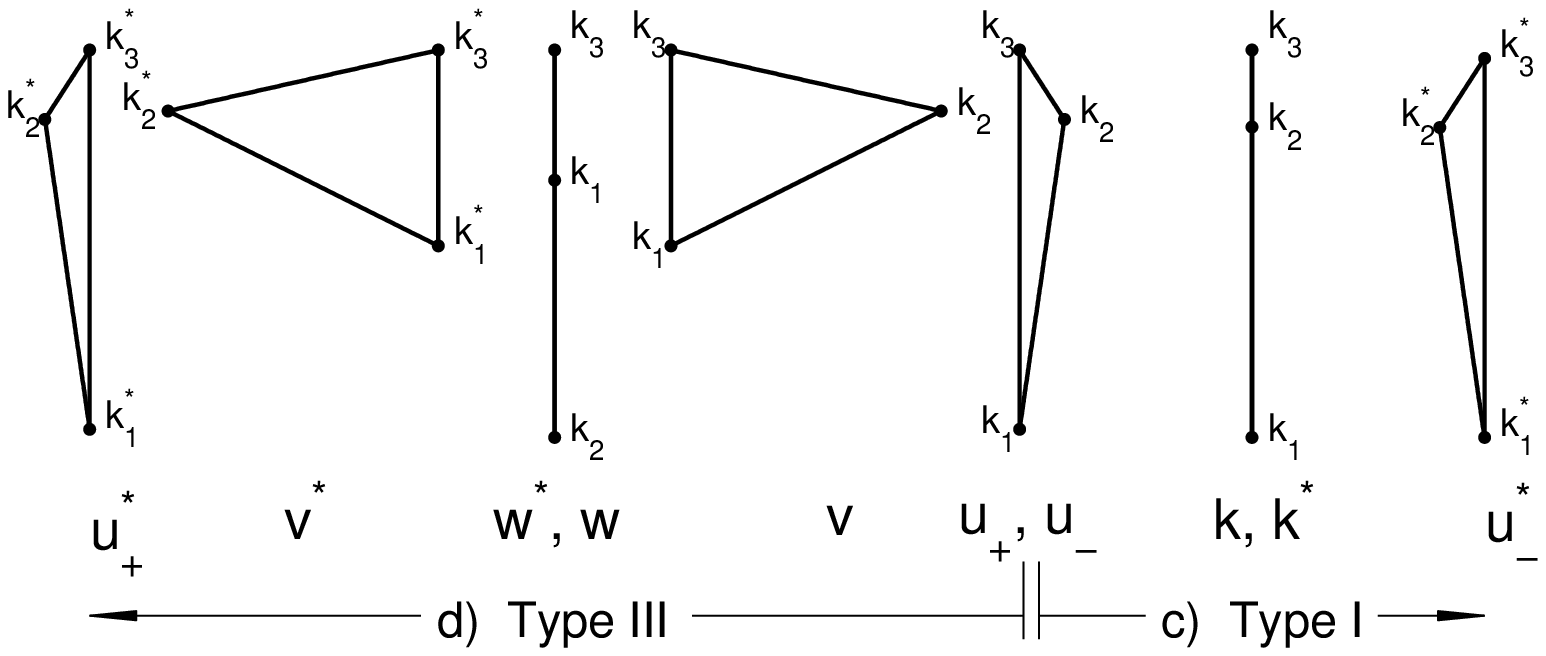}}
\caption{Evolution of $\triangle_1$ departing 
from the unstable branch
on the positive side of $\cal P$ to the stable branch on the negative
side of type I or III, or to that on the positive side of type II.}
\label{fig:shapeAB}
\end{figure}

Figures~\ref{fig:twoby2-Y}, \ref{fig:twoby2-p}
and \ref{fig:shapeAB} show, respectively, the variations of 
the perimeter, $\log p(t)$, the deviation from the critical curve,
$Y(t)$, and the vortex configuration 
scaled by its perimeter, $\triangle_1$, along the trajectories 
for the four cases listed in  Table~\ref{table:1}.

When a trajectory crosses over an edge to the opposite side of $\cP$,
i.~e., $\gamma$ changes from $\pm 1$ to $\mp 1$, it continues along
the image of that on the positive side but in reverse direction;
therefore,  at an edge,  $\log p$  reaches  an extremum and so does $Y$
while the three vortices become collinear.
We also note that   on either side of
$\cP$ or the $\alpha \beta$ plane, 
the trajectories,
$\bar I =$ constant, are concave
downward while the constant $Y(x_j)$ lines are concave upward
similar to that of $\cC$ where $Y=0$.
The trajectory  $\bar I = 1$ and the critical curve $Y=0$ are 
of opposite concavity. They are tangent
to each other a $E$ (or $E^*$).
The trajectories intersecting the left (right) 
branch $Q_4 E$ or $EQ_5$ of $\cC$ will have
$\bar I(Q_4) < \bar I < 1$ or $\bar I(Q_5) < \bar I < 1$. This statement holds
also for the images of the trajectories.

The trajectories with $\bar I < 1$
will be tangent to a constant $Y$-line with $Y > 0$ above $\cC$. Thus
along  a trajectory with $\bar I < 1$, $Y(\bx)$ reaches
a maximum at a point above $\cC$. Again this holds also for the
image trajectory on the negative side of $\cP$

Along a  trajectory initiating  at point
below $\cC$, say $r^-$ or $u^-$, the deviation,  $\delta = - Y$, 
(the perimeter
or $- \log p$) will increase from 
almost zero, or $\delta \ll 1$,   to a maximum  at the edge and then
return to zero as the trajectory  ending at a point 
on the stable branch of $\cC$ while $\triangle_1 (t)$ becomes equivalent
to $\triangle_1(0)$ with a change of orientation. 
The trajectory is of type I. See 
part $a)$ and $c)$ of these three figures.

For a trajectory initiating from a point, say $r^+$,  above the unstable branch 
$ES_4$ of $\cC$ where $\bar I(S_4) = \bar I (Q_4) >  \bar I (Q_5)$,
the deviation  $\delta$   will
increase from a small initial value , $\ll 1$, to a maximum and
then return to zero as 
the trajectory  arrives at a point on
the stable branch $Q_4E$  of $\cC$ without crossing over
an edge, i.~e., staying on the same side of $\cP$, as  
in  Fig.~\ref{fig:twoby2-Y} b). The trajectory 
is of type II. The perimeter will decrease from $1$
to a minimum and then keep on increasing as $t$ increases.
See part b) in Fig.~\ref{fig:twoby2-p}.
Note that  the final configuration corresponding to the point
on the left branch $Q_4E$
is different from the initial one on the right branch $ES_4$, as shown
in Fig.~\ref{fig:shapeAB} d).

For a trajectory initiated from a point, $u^+$ above $\cC$ with $\bar I \in
(\bar I (Q_5), \bar I (Q_4))$,  the 
deviation, $\delta = Y(0)$,  will
increase from a  small initial value $ \ll 1$ to a maximum, crossing
the edge $Q_2Q_3$ with $Y$ reaching a local minumum and continue
to the negative side of $\cP$ with $Y$ appearing 
as the image of the first part passing through the image point of $u^+$, 
then return to zero as
it reaches a point on
the stable branch of $\cC$. The image of the point will be
close to the initial point $u^+$. The trajectory
is of type III.

Our numerical examples show the three distinct  types of
trajectories expected in the theoretical studies in the preceding
section.
Note that the deviation $Y(t)$  and the perimeter $p(t)$
for type III have three extrema as $t$ increases while
those for type I and II have only one.
For type I and III, there are changes of orientation  
and the final configurations are nearly similar to the original ones
with perimeter  increasing. For type II, there
are no changes of orientation but the final expanding configuration
is not similar to the initial one.

The trajectories depart from the unstable branch
$E^*Q_4$ of $\cC$ on the negative side of $\cP$. Because 
of $\bar I_4 > \bar I_5$, the departing trajectories
are only  of two types: trajectories of
type I begin from a point below the branch $Q_4E^*$
while those of type II begin from a point above.

\section{Conclusion}
\label{sec-conclusion}                 

\noindent 
This paper presents a complete description of
the grobal dynamics of three point vortices initially in the neighborhood of 
the planar critical curve $\cC$ in trilinear coordinates $x_j, \ j=1, 2, 3$.
The trilinear coordinates, only two of which  are independent because of
the constraint $x_1 + x_2 + x_3 = 1$,  were introduced
by Synge (1949) to study the deformation of the triangle
$\triangle_p$ formed
by the three vortices with sides $R_j$ and perimeter $p$
via that  of a similar triangle with sides $x_j$ and perimeter one.
A spatial trajectory $R_j(t)$ is then reduced to a planar trajectory $x_j(t)$.
In the plane, a stationary critical curve $\cC$
was found by Synge (1949) for parabolic cases.
Each point on $\cC$ corresponds to a radial spatial trajectory
either moving away or towards the  origin with similar $\triangle_p$
expanding $\dot p > 0$, or contracting 
$\dot p < 0$, respectively. The latter would lead to 
coelescence of three vortices. According to the spatial trajectories,
the  critical curve is 
partitioned into expanding or contracting branches. It was shown by Tavantzis 
and Ting (1988) that the
expanding radial (spatial)  trajectories are stable while
the contracting ones are unstable. In this paper we show that
a  contracting radial trajectory slightly disturbed will depart from the
radial line,  move along a spatial trajectory
nonsimilar to the initial configuration and eventually approach
a radial expanding trajectory. Thus we come to the conclusion
that three distinct vortices almost never coelescence.


\begin{thebibliography}{11}
\bibitem{grobli}  W. Gr\"{o}bli, Specielle Prcbleme \"{u}ber die Bewegung
geradliniger paraller Wirbelf\"{a}den. Z\"{u}rich: Z\"{u}rich and Furrer,
1877.

\bibitem{synge} J. Synge,  On the motion of three vortices. Can. J. Math. 1949;
1; 257-270.

\bibitem{novikov}  E. Novikov, Dynamics and statistics of a system of vortices. 
Sov.  Phys. JETP 1975; 41; 937-943.

\bibitem{aref}  H. Aref, Motion of three vortices. Phys. Fluids 1979; 22;
393-400.

\bibitem{tavtin}  J. Tavantzis, L. Ting, The dynamics of three vortices revisited.
Phys. Fluids 1988; 31; 1392-1409.

\bibitem{knicol}  O. Knio, L. Collorec, D. Juv\'{e}, Numerical 
study of sound emissions by 2D regular and chaotic configurations. 
J. Comput. Phys. 1995; 116; 226-246.

\bibitem{BTK06} D. Blackmore, L. Ting, O. Knio,
``Studies of perturbed three vortex dynamics'', J. Math. Phys.,
Vol. 48, 065402, pp. 1--30, 2007.

\bibitem{TKB07} L. Ting, O. Knio, D. Blackmore,
``Dynamics of Planar Vortex Clusters  with  Binaries''
presented in Minisymposium IC/MP/107/R/520, 2nd part, 
``Recent Advances in Vortex Dynamics: Theory and Computation'',
ICIAM-07, Zurich, July 16-20, 2007. To appear in
the Proc. ICIAM-07.

\bibitem{KBT08} Knio, D. Blackmore, L. Ting,
``Numerical Study of Dynamics of Point Vortex Configurations'',
presented at  ICCES08, 17-22 March 2008, Hawaii.

\bibitem{lin43} C. C. Lin, ``On the motion of vortices 
in two-dimensions'', Toronto Univ. Press, Toronto, 1943.


\bibitem{newton} P. Newton, The N-Vortex Problem, Springer, 
New York, 2001.

\bibitem{tinbla}  L. Ting, D. Blackmore, Bifurcation of motions of three vortices and applications,  presented
in ICTAM-04, Session FM25, {\it Vortex Dynamics},
Warsaw, Poland, XXI ICTAM 2004 Abstracts and  CD-ROM Proceedings, 
pp. 188-189.


\bibitem{lamb} H. Lamb, Hydrodynamics, 6th ed., Dover republ. New York, 1945.

\end{thebibliography}
\end{document}